\documentclass[a4paper]{article}
\usepackage{amsthm,amsmath,amssymb}
\usepackage{color}
\newtheorem{teor}{Theorem}[section]
\newtheorem{defin}[teor]{Definition}
\newtheorem{lemm}[teor]{Lemma}
\newtheorem{osse}[teor]{Remark}
\newtheorem{prop}[teor]{Proposition}
\newtheorem{defi}[teor]{Definition}
\newtheorem{coro}[teor]{Corollary}
\newtheorem{prob}[teor]{Problem}

\newcommand{\bele}{\begin{lemm}\begin{sl}}
\newcommand{\enle}{\end{sl}\end{lemm}}
\newcommand{\bedef}{\begin{defi}\begin{sl}}
\newcommand{\eddef}{\end{sl}\end{defi}}
\newcommand{\bete}{\begin{teor}\begin{sl}}
\newcommand{\ente}{\end{sl}\end{teor}}
\newcommand{\beos}{\begin{osse}\begin{rm}}
\newcommand{\eddos}{\end{rm}\end{osse}}
\newcommand{\bepr}{\begin{prop}\begin{sl}}
\newcommand{\empr}{\end{sl}\end{prop}}
\newcommand{\bepro}{\begin{prob}\begin{rm}}
\newcommand{\empro}{\end{rm}\end{prob}}
\newcommand{\bede}{\begin{defin}\begin{sl}}
\newcommand{\edde}{\end{sl}\end{defin}}
\newcommand{\beco}{\begin{coro}\begin{sl}}
\newcommand{\enco}{\end{sl}\end{coro}}

\newcommand{\quext}{\quad\text}
\newcommand{\qquext}{\qquad\text}
\newcommand{\de}{\partial}

\definecolor{grey}{rgb}{0.85,0.85,0.85}

\long\def\greybox#1{%
    \newbox\contentbox%
    \newbox\bkgdbox%
    \setbox\contentbox\hbox to \hsize{%
        \vtop{
            \kern\columnsep
            \hbox to \hsize{%
                \kern\columnsep%
                \advance\hsize by -2\columnsep%
                \setlength{\textwidth}{\hsize}%
                \vbox{
                    \parskip=\baselineskip
                    \parindent=0bp
                    #1
                }%
                \kern\columnsep%
            }%
            \kern\columnsep%
        }%
    }%
    \setbox\bkgdbox\vbox{
        \color{grey}
        \hrule width  \wd\contentbox %
               height \ht\contentbox %
               depth  \dp\contentbox
        \color{black}
    }%
    \wd\bkgdbox=0bp%
    \vbox{\hbox to \hsize{\box\bkgdbox\box\contentbox}}%
    \vskip\baselineskip%
}

\newcommand{\RR}{\mathbb{R}}

\newcommand{\beeq}[1]{\begin{equation}\label{#1}}
\newcommand{\eddeq}{\end{equation}}

\newcommand{\beeqa}[1]{\begin{eqnarray}\label{#1}}
\newcommand{\eddeqa}{\end{eqnarray}}

\newcommand{\beal}[1]{\begin{align}\label{#1}}
\newcommand{\eddal}{\end{align}}

\newcommand{\bespl}[1]{\begin{split}\label{#1}}
\newcommand{\edspl}{\end{split}}

\newcommand{\bega}[1]{\begin{gather}\label{#1}}
\newcommand{\edga}{\end{gather}}

\newcommand{\beeqax}{\begin{eqnarray*}}
\newcommand{\eddeqax}{\end{eqnarray*}}

\def\qed{\ifmmode 
  \else \leavevmode\unskip\penalty9999 \hbox{}\nobreak\hfill
  \fi
  \quad\hbox{\hskip.5em\vrule width.4em height.6em depth.05em\hskip.1em}}
\def\endproofsym{\qed}
\renewenvironment{proof}[1][Proof]{\trivlist\item[\hskip\labelsep{\hskip0pt
    {\normalfont\scshape#1.}\hskip .321429\parindent}]\ignorespaces}
{\endproofsym\endtrivlist}
\def\endnobox{\def\endproofsym{}\end{proof}\def\endproofsym{\qed}}

\newcommand{\no}{\nonumber}

\newcommand{\beeqao}{\begin{eqnarray}\no}
\newcommand{\bealo}{\begin{align}\no}
\newcommand{\besplo}{\begin{split}\no}
\newcommand{\begao}{\begin{gather}\no}

\newcommand{\eps}{\epsilon}

\newcommand{\duav}[1]{\langle{#1}\rangle}

\newcommand{\duavw}[1]{\langle\!\langle{#1}\rangle\!\rangle}

\newcommand{\ov}[1]{\overline{#1}}

\newcommand{\perogni}{\forall\,}
\newcommand{\esiste}{\exists\,}
\newcommand{\itt}{\int_0^t}

\newcommand{\io}{\int_\Omega}

\newcommand{\ito}{\itt\!\io}

\newcommand{\iTT}{\int_0^T}

\newcommand{\iTo}{\iTT\!\io}

\newcommand{\epsi}{\varepsilon}
\newcommand{\ee}{^{\varepsilon}}

\newcommand{\zt}{_{(0,t)}}

\newcommand{\tilfhi}{\widetilde{\fhi}}

\def\R{\mathbb R}

\newcommand{\bn}{\boldsymbol{n}}

\newcommand{\fhi}{\varphi}

\newcommand{\lhs}{left hand side}
\newcommand{\rhs}{right hand side}

\newcommand{\devv}{\partial_{{\calV},{\calV'}}}

\DeclareMathOperator{\deriv}{d}

\newcommand{\LDH}{L^2(0,T;H)}

\let\TeXchi\chi
\def\chi{{\setbox0 \hbox{\mathsurround0pt
$\TeXchi$}\hbox{\raise\dp0 \copy0 }}}

\newcommand{\calH}{{\mathcal H}}

\newcommand{\calX}{{\mathcal X}}

\newcommand{\calT}{{\mathcal T}}

\newcommand{\calE}{{\mathcal E}}
\newcommand{\calJ}{{\mathcal J}}

\newcommand{\calV}{{\mathcal V}}
\newcommand{\Vtzero}{{\mathcal V}_{t,0}}

\newcommand{\Xtzero}{{\mathcal X}_{t,0}}

\newcommand{\dit}{\deriv\!t}
\newcommand{\dis}{\deriv\!s}
\newcommand{\dix}{\deriv\!x}

\newcommand{\diT}{\deriv\!{\calT}}
\newcommand{\dir}{\deriv\!r}

\newcommand{\ddt}{\frac{\deriv\!{}}{\dit}}

\newcommand{\ditau}{\deriv\!\tau}

\newcommand{\bFormula}[1]{
\begin{equation} \label{#1}}
\newcommand{\eF}{\end{equation}}

\newcommand{\calD}{{\mathcal D}}

\definecolor{violet}{rgb}{0.85,0.05,0.85}

\newenvironment{betti}{\color{red}}{\color{black}}
\newcommand{\bebe}{\begin{betti}}
\newcommand{\ebe}{\end{betti}}

\newenvironment{ele}{\color{blue}}{\color{black}}
\newcommand{\elena}{\begin{ele}}
\newcommand{\finele}{\end{ele}}

\newenvironment{riccardo}{\color{green}}{\color{black}}
\newcommand{\beri}{\begin{riccardo}}
\newcommand{\eri}{\end{riccardo}}

\newenvironment{giu}{\color{violet}}{\color{black}}
\newcommand{\giulio}{\begin{giu}}
\newcommand{\fingiu}{\end{giu}}

\numberwithin{equation}{section}
\begin{document}

\title{On the strongly damped wave equation with 
 constraint}

\author{%
Elena Bonetti\\
Dipartimento di Matematica, Universit\`a di Pavia,\\
Via Ferrata~1, 27100 Pavia, Italy\\
E-mail: {\tt elena.bonetti@unipv.it}
\and
Elisabetta Rocca\\
Weierstrass Institute for Applied
Analysis and Stochastics\\ Mohrenstr.~39, 10117 Berlin, Germany\\ 
E-mail: {\tt elisabetta.rocca@wias-berlin.de} \\
and \\
Dipartimento di Matematica, Universit\`a di Milano,\\
Via Sal\-di\-ni 50, 20133 Milano, Italy\\
E-mail: {\tt elisabetta.rocca@unimi.it} 
\and 
Riccardo Scala\\
Dipartimento di Matematica, Universit\`a di Pavia,\\
Via Ferrata~1, 27100 Pavia, Italy\\
E-mail: {\tt riccardo.scala@unipv.it}
\and 
Giulio Schimperna\\
Dipartimento di Matematica, Universit\`a di Pavia,\\
Via Ferrata~1, 27100 Pavia, Italy\\
E-mail: {\tt giusch04@unipv.it}
}


\maketitle
\begin{abstract}
 A weak formulation for the so-called {\sl semilinear strongly 
 damped wave equation with constraint}\/ is introduced and a corresponding
 notion of solution is defined. The main idea in this approach consists in
 the use of duality techniques in Sobolev-Bochner spaces, aimed at
 providing a suitable ``relaxation'' of the constraint term. 
 A global in time existence result is proved 
 under the natural condition that the initial data have 
 finite ``physical'' energy.
\end{abstract}

\noindent {\bf Key words:}~~wave equation, strong damping, weak
 solution, maximal monotone operator, duality.

\vspace{2mm}

\noindent {\bf AMS (MOS) subject clas\-si\-fi\-ca\-tion:} 
35L05, 74D10, 47H05, 46A20.


\section{Introduction}
\label{sec:intro}

This paper is devoted to studying the so-called semilinear
wave equation with strong damping, namely
\begin{equation}\label{eqn:intro}
  \epsi u_{tt} 
   - \delta \Delta u_t
   - \Delta u
   + f(u) 
   = g,
\end{equation}
for $\epsi,\delta >0$. The equation is
settled in the parabolic cylinder 
$Q=(0,T)\times \Omega$, where $\Omega$ is a smooth
bounded domain in $\RR^N$, $N\ge 1$, and $T>0$
is a given final time, and is complemented
with the initial conditions for $u$ and $u_t$ and
with homogeneous boundary conditions either of 
Dirichlet or of Neumann type.
The {\sl strong damping}\/ is 
provided by the term $- \delta \Delta u_t$;
this comes in contrast with the {\sl weak damping}\/
occurring when that term is replaced by $+ \delta u_t$.
The function $g$ on the \rhs\ is a given volume forcing term (here taken  
of $L^2$-regularity), and the semilinear term 
$f(u)$ is assumed to take the form $f(u) = \beta(u) - \lambda u$, where 
$\beta$ is a monotone function (more precisely,
a monotone {\sl graph}, see Section~\ref{sec:main}
below) and $\lambda \ge 0$. In particular, the
internal constraint on $u$ is enforced by the non-smooth
monotone part $\beta$ of $f$, whereas the remaining
term $- \lambda u$ is related to the (possible)
nonconvexity of the energy functional
associated to the equation.
Actually, the main novelty of this paper stands in the fact
that $\beta$ is assumed to be defined only
in a {\sl bounded}\/ interval $I_0$ of $\RR$ and 
to diverge at the extrema of $I_0$. A 
(generalized) function $\beta$ with the 
above properties will be referred to as 
a {\sl constraint}\/
on the variable $u$ (cf.~Section~\ref{sec:main}
below for more details). It is worth noting that,
up to purely technical modifications in the proofs,
our techniques could be adapted to treat also the case of
{\sl unilateral constraints}, i.e., functions 
$\beta$ whose domain is bounded only from one side. 

Physically speaking, equation~\eqref{eqn:intro} appears in a number of 
different contexts. Let us mention here some of them. 
The main application refers to
the study of the motion of viscoelastic materials. 
In this setting, $u$ plays the role of a (scalar) displacement and \eqref{eqn:intro}
represents the momentum balance (where accelerations are included) written in a  small strain regime.  
In particular, respectively in space dimensions 
one and two, the equation describes
the transversal vibrations of a homogeneous string and the longitudinal 
vibrations of a homogeneous bar subject to viscous effects. 
The strong damping term $- \delta \Delta u_t$ represents the fact that the stress 
is decomposed in the sum of a pure elastic part (proportional to the strain) and a viscous 
part (proportional to the strain rate), as in a linearized Kelvin-Voigt material.
We also mention that in the literature, in space dimension three, 
\eqref{eqn:intro} has been introduced to model, e.g., 
the deviation from the equilibrium configuration of a 
(homogeneous and isotropic) linearly viscoelastic 
solid with short ``rate type'' memory (cf.~\cite{DL} for details), 
in the presence of an external displacement-dependent force $g-f(u)$.
We do not enter deeper in the modeling details, and we 
refer to \cite{K} for a physical derivation of models describing
the motion of viscoelastic media. 
Let us observe that it would be meaningful to consider
here a vectorial (displacement) variable
${\bf u}$, but we preferred, just for simplicity, 
to study only the scalar case at least at a first stage.
Indeed, the extension of our results to the vector-valued case should 
be possible, at least for constant isotropic
diffusion, whereas the case of non-constant stiffness  
(and viscosity) tensors may be somehow more involved. In this
framework, we also have to quote (possibly adhesive) contact 
models with unilateral constraints 
(occurring for instance in the case of Signorini conditions) 
on a part of the boundary. In this setting, the
(vectorial) operator $\beta$ would
force the direction of the trace of ${\bf u}$ on the boundary 
in such a way to ensure impenetrability (cf., e.g., 
\cite{BBR1,BBR2,RR}). We are planning to analyze this type
of models, by using the methods developed in this paper,
in future works.

Equation \eqref{eqn:intro} also appears in the so-called Fr\'emond theory
for phase transitions whenever microscopic accelerations are taken into account
(cf., e.g.,~\cite{bl1,bl2,FR}). In that setting, the unknown $u$ generally 
denotes a (scalar) phase parameter, which is related (for a first 
order phase transition in a binary system) to
the local proportion of one of the two phases, or components,
of a binary material. Then, $\beta$ 
represents an internal constraint forcing 
$u$ to take values into the physical interval 
whose extrema (often given by $-1,1$) correspond to the pure 
states, whereas the intermediate values represent a mixture
of the phases. Physically relevant choices are 
$\beta(r)=\log (1+r) - \log(1-r)$ (i.e., the derivative of the 
so-called {\sl logarithmic potential}\/ often appearing in Allen-Cahn
or Cahn-Hilliard models), or $\beta(r)=\de I_{[-1,1]}(r)$
(i.e., the {\sl subdifferential}\/ of the 
{\sl indicator function}\/ of the interval $[-1,1]$, given by
$I_{[-1,1]}(r)=0$ for $r\in [-1,1]$ and
$I_{[-1,1]}(r)=+\infty$ otherwise). 
It is also significant to consider equation \eqref{eqn:intro} 
with other kinds of nonlinearities $\beta$ not having the form of 
a constraint. For instance, \eqref{eqn:intro} 
appears in the recent theory of isothermal 
viscoelasticity with very rapidly fading memory (cf.~\cite{CP} 
and references therein), in the sine-Gordon model describing the 
evolution of the current $u$ in a Josephson junction 
(cf.~\cite{LSC}; there $f(u)=\sin u$), or as a 
Klein-Gordon-type equation occurring in quantum mechanics 
(then $f(u)=|u|^\gamma u$ for suitable $\gamma > 0$).
%
%

Actually, in the case when $f$ is smooth and defined on the whole real line,
the mathematical literature on equation \eqref{eqn:intro}
is very wide (we quote, without any claim of completeness, the papers
\cite{AS,GM,HALE,KAL,KZ,LE,PZ,WEBB}). Referring to \cite{KZ} for
more details, we recall here that one of the first essential 
results on global well-posedness of (the Dirichlet problem for) 
\eqref{eqn:intro} in the 3D case was obtained by Webb,
who proved in \cite{WEBB} that, if $f$ satisfies standard
dissipativity conditions (without any growth restriction), 
then the problem admits a unique strong solution $(u, u_t)$ 
taking values in the space $[H^2(\Omega)\cap H_0^1(\Omega)] \times L^2(\Omega)$. 
On the other hand, when one looks for less regular solutions, the situation
seems different. In particular, it is natural to consider weaker
solutions such that the ``energy of the system'' remains bounded 
(in the analytical literature this fact corresponds
to require that these solutions take values in the so-called 
{\sl energy space}\/). Indeed, this type of regularity corresponds to the a priori estimate
obtained by (formally) testing \eqref{eqn:intro} by $u_t$. 
Then, one can easily realize that, at least if 
the external source $g$ is $0$, the functional
\begin{equation}\label{ener:intro}
  \calE(u,u_t) 
   = \io \Big(
    \frac{\epsi}2 | u_t |^2 
    + \frac12 | \nabla u |^2
    + j(u) - \frac\lambda2 u^2 \Big)\,\dix,
\end{equation}
where $j$ is an antiderivative of $\beta$, tends to decrease in the time 
evolution. Usually $\calE$ is interpreted as a physical energy. This is 
particularly clear in the cases when $u$ represents a displacement
(including phase-change models where $u$ is related to the effects of
displacements at microscopic scales): then
the component $|u_t|^2$ of the integrand is a density of kinetic energy, 
whereas the other summands correspond to some kind of configurational
or potential energy. 
Consequently, energy solutions can be defined as those solutions 
taking values in the energy space, or, equivalently, 
keeping finiteness of the energy in the course of the evolution.

From the mathematical point of view, managing this type of solutions 
may be delicate, especially in high space dimension, 
in view of the possibly fast growth of the integrand $j(u)$. 
Correspondingly, the literature related
to this case is much more recent: Kalantarov and Zelik in \cite{KZ}
consider polynomial nonlinearities of the form 
$f (u)\sim  u|u|^q$ without any restriction on the exponent $q>0$
and prove well-posedness of the equation in the energy space.
Moreover, they analyze the asymptotic behavior of solutions for
large values of the time variable. More recently Pata and Zelik in~\cite{PZ} 
have extended these results to the 
case when $f:\RR\to \RR$ is any smooth function satisfying the 
basic dissipativity assumptions, without any restriction on the growth
rate (for instance, they may take exponential, or even faster growing,
terms $f$). At least up to our knowledge, however, the case when 
$f$ is of {\sl constraint}\/ type has never been considered up to now. 
The typical example we have in mind is $f(u)=\partial I_{[-1,1]}(u) - u$ 
(cf.~\eqref{hp:f} and \eqref{hp:beta} below), which may describe,
for example, some phase transition phenomena accounting for
microscopic accelerations (cf.~\cite{FR}). 

As one addresses the initial-boundary value problem for~\eqref{eqn:intro}
under our assumptions, the main mathematical difficulty   
comes from the combination 
of the constraint $\beta$ with the 
second time derivative~$u_{tt}$. Indeed, 
this feature strongly restricts the available a-priori 
bounds. To be precise, almost all information
on the solution has to be extracted from the so-called  
``energy'' estimate, i.e., testing the equation by $u_t$. 
In addition to that, one can just get some more smoothness
of $u$ by multiplying \eqref{eqn:intro} by $-\Delta u$
(as is done, e.g., in~\cite{PZ}). Anyway, this 
does not help for controlling the term $u_{tt}$, which is the main
issue from the point of view of regularity. Moreover, 
the standard procedures that one usually adopts for 
obtaining higher order bounds, like differentiating in time the 
equation, do not seem to work here, at least for a general choice of $\beta$. 
This seems to be, indeed, the main difference of the present problem with
respect to
{\sl first order}\/ (in time) 
equations with constraint,
for which additional regularity of solutions can be generally 
deduced by differentiating in time and testing 
the result by $u_t$, whatever is the expression of $\beta$.

In view of the lack of estimates, we need to build a notion
of weak solution which is sufficiently general to exist under
the sole ``energy'' regularity.
This is, indeed, a somehow delicate issue.
In particular, one cannot expect to reproduce the same type of 
results that hold in the case of less
general nonlinearities $\beta$. To say it shortly,
the main novelties of our approach can be summarized in 
two points:\\[2mm] 
$(i)$~~a relaxed form $\beta_w$ of the operator $\beta$ obtained by means of 
 duality techniques;\\[1mm]
$(ii)$~~an integrated (both in space and in time)
 variational formulation where test functions are chosen in suitable 
 Sobolev-Bochner spaces.\\[2mm]
These choices permit us, indeed, to prove existence. However, 
both of them come at some price. Namely,\\[2mm]
$(i)$~~it will not be possible to intend the equation, and the 
constraint in particular, in the pointwise sense;\\[1mm]
$(ii)$~~we cannot exclude the occurrence of {\sl jumps}\/ of $u_t$.
 Actually, $u_t$ may be discontinuous with respect to time (and, more precisely,
 is {\sl expected}\/ to be discontinuous, as we can show by means of 
 examples).\\[2mm]
However, from a physical point of view, if \eqref{eqn:intro} 
comes from a  variational principle (as the principle of virtual power is), 
the variational setting in which we introduce the solution
is the natural one. In particular, the operator $\beta$ (in its
relaxed version $\beta_w$) stands for an internal force which is defined
in duality with velocities/displacements. In addition to that, 
an internal constraint on the function~$u$ is still ensured by
the definition of the domain of $\beta_w$. Finally, the 
fact that we can have jumps on the velocity $u_t$ w.r.t. time, corresponds
to the possible occurrence of internal (or external) shocks, which are
expected to happen in this framework 
(cf, e.g., \cite{FR}). 

A further drawback is concerned with the problem of uniqueness.
Actually, we expect the occurrence of genuine nonuniqueness, even
though some criteria for ``physicality'' of weak solutions may be 
proposed (cf.~Remark~\ref{rem:reg2} at the end). 

Let us conclude by giving some more words of explanation for our method.
The basic strategy of proof is, in a sense, very standard: we replace the 
singular function $\beta$ by a smooth approximation $\beta\ee$ of controlled
growth at infinity (e.g., the Yosida approximation), prove existence
of a solution $u\ee$ to the regularized problem
(which basically follows from results already known in
the literature, cf., e.g., the quoted \cite{GM,KZ,PZ,WEBB}),
and then let the approximation parameter $\epsi$ go to $0$.
Indeed, as a consequence of the so-called ``energy estimate'',
$u\ee$, at least for a subsequence, tends to some limit $u$
which we would like to identify as a ``weak'' solution
(where, of course, we need to state precisely what we 
mean with this).
However, the only uniform bound available for the nonlinear term 
$\beta\ee(u\ee)$ is in the
norm of $L^1(0,T;L^1(\Omega))$, and there is evidence coming
from concrete examples that we cannot go further,
at least for general $\beta$. This fact has, indeed, a number of consequences. 
First of all, arguing by comparison, 
we can obtain an $L^1(0,T;X)$-bound for $u\ee_{tt}$, 
where $X$ is a Banach space such that $L^1(\Omega)$ and $H^1(\Omega)$
are compactly embedded into $X$ (for example,
we can take $X=H^{-2}(\Omega)$ in the 3D case).
This estimate suffices, via a generalized version of the Aubin-Lions
compactness lemma, to prove strong convergence of $u\ee_t$ 
in $L^2(0,T;L^2(\Omega))$. However, the limit function $u_t$ 
may exhibit jumps with respect to time. Secondly, the limit
of $\beta\ee(u\ee)$ can be taken at least in the (weak)
sense of measures. A crucial point is, as usual, concerned
with the identification of its limit. In view of our
assumptions it looks natural to rely on a suitable version of the
so-called Minty's trick for monotone operators,
i.e., to combine the weak convergence of $u\ee$ in some (reflexive)
Banach space $\calV$, the weak convergence of $\beta\ee(u\ee)$ 
in the dual space $\calV'$, and a $\limsup$-inequality. 
A look at the estimates suggests that an admissible
choice for this procedure is (in the Dirichlet case)
the Sobolev-Bochner space $\mathcal{V}=H^1(0,T;L^2(\Omega))\cap L^2(0,T;H^1_0(\Omega))$
(in the Neumann case, $H^1_0(\Omega)$ is simply replaced 
by $H^1(\Omega)$ and no further difficulties arise).
In such a setting, the constraint $\beta$ has to be 
reinterpreted in a relaxed form $\beta_w$ acting as a 
maximal monotone operator from $\calV$ to $2^{\mathcal{V}'}$ 
(cf.~Definition~\ref{defibw} below; see, e.g., \cite{brezisart,GR}
for some additional background). 
Correspondingly, equation~\eqref{eqn:intro} has to be intended as a relation in $\mathcal{V}'$. 
Let us point out that, from a physical point of view, in the case
when \eqref{eqn:intro} corresponds to a mechanical balance equation (i.e.,
to the momentum balance equation),
our weak formulation takes the meaning of a duality between 
forces and velocities in time and space 
(see \cite{BCFlarge} for a similar approach, but in a different setting). 
To avoid occurrence of second time derivatives in the weak formulation,
we also need to integrate by parts with respect to time the second order
term $u_{tt}$ (cf.~\eqref{eqnw}). Actually, these modifications
will permit us to solve our original problem on the whole time interval $(0,T)$, 
but also to write ``pointwise'' the duality 
relation in any subinterval $(0,t)$, with a physically 
consistent interpretation of the corresponding constraint. 
Finally, an energy inequality is proved
to hold on (almost) every subinterval of~$[0,T]$. 
We end observing that the behavior of weak solutions
(at least in the homogeneous Neumann case)
may be clarified by considering a spatially homogeneous
setting. For instance, if $f(u)=\partial I_{[-1,1]}(u)$ and $g\equiv 0$, 
\eqref{eqn:intro} reduces to the prototype ODE
$u_{tt}+\partial I_{[-1,1]}(u)\ni 0$ whose solutions
can be easily described, especially in relation with the jumps
of $u_t$ (cf.~Remark~\ref{rem:reg1} for more details). 

\smallskip

The remainder of the paper is organized as follows: 
in Section~\ref{sec:main} we introduce some amount of 
preliminary material mainly related to maximal monotone operators
and duality methods; moreover we present the
notion of weak solution and state the related existence
result. Then, the proof is detailed in 
Section~\ref{sec:proofs}, where we also give a number
of remarks illustrating our results at the light of 
simple finite-dimensional examples.


\section{Preliminary notions and main result}
\label{sec:main}

Let $\Omega\subset\RR^N$ be a smooth bounded domain (with $N\geq1$)
of boundary~$\Gamma$ and let us consider the interval $[0,T]$, 
for some fixed final time $T>0$. Let us set $H:=L^2(\Omega)$ 
and use the notation $(\cdot,\cdot)$ 
for the scalar product both in $H$ and in $H^N$. Let also the symbol
$\| \cdot \|$ denote the corresponding
norms. In our analysis, we will consider either Dirichlet or Neumann boundary
conditions for \eqref{eqn:intro}; hence we introduce a notation suitable 
for addressing both cases in a unified way. So, we put
$V:=H^1(\Omega)$ in the Neumann case, and $V:=H^1_0(\Omega)$ in the Dirichlet
case.
%
%
%
%
In both cases, $V$ will be endowed with the standard (Sobolev) norm,
indicated by~$\| \cdot \|_V$. Moreover, 
we will denote by $\duav{\cdot,\cdot}$
the duality pairing between $V'$ and $V$.
In general, we will indicate by $\| \cdot \|_X$
the norm in some Banach space $X$
(or in $X^N$).
%

We let $A$ stand for the weak form of (minus) the Laplace operator
seen as an unbounded linear operator
on $H$ whose domain $D(A)$ depends on the boundary conditions.
Namely, in the Neumann case, we set
\begin{equation}\label{Aneum}
  Av := -\Delta v,
   \quad D(A):= H^2_{\bn}(\Omega),
\end{equation}
where $H^2_{\bn}(\Omega)$ denotes the space
of the $H^2$-functions having zero normal derivative (in the 
sense of traces) on $\de\Omega$. 
Correspondingly, in the Dirichlet case, we set
\begin{equation}\label{Adir}
  Av := -\Delta v,
   \quad D(A):= H^2(\Omega) \cap V
  = H^2(\Omega) \cap H^1_0(\Omega).
\end{equation}
In both cases, $A$ is a positive operator (strictly
positive for Dirichlet conditions). Morever, $A$
can be extended to the space $V$ by setting
(for both choices of boundary conditions)
\begin{equation}\label{Aduav}
  \duav{A v,z} = \io \nabla v\cdot\nabla z\,\dix. 
\end{equation}
This extension, which turns out to be 
linear and bounded from $V$ to $V'$, will 
be identically noted as $A$; indeed, we believe
that no danger of confusion exists at this stage. 

Next, we specify our assumptions on the semilinear term
$f(u)$. First, we suppose that $f$ may be decomposed as
\begin{equation}\label{hp:f}
  f(u) = \beta(u) - \lambda u,
\end{equation}
where $\lambda \ge 0$ and $\beta$ is a 
{\sl maximal monotone graph}\/ 
in $\RR\times \RR$ such that 
\begin{equation}\label{hp:beta}
  \ov{D(\beta)}=[-1,1], \quad 0\in \beta(0).
\end{equation}
Indeed, just for simplicity and with
no loss of generality, we require the 
closure of the {\sl domain}\/ of $\beta$ to be 
the interval $[-1,1]$. In addition, it is not restrictive to assume 
the normalization $0\in \beta(0)$, which turns out to be
useful especially in the Dirichlet case.

Referring the reader to \cite{Ba,Br} for a complete
survey on the theory of maximal monotone operators
in Banach and Hilbert spaces, we just observe here that,
thanks to \eqref{hp:beta}, there exists a convex
and lower semicontinuous function 
$j:\RR\to [0,+\infty]$ such that $\beta=\de j$,
$\ov{D(j)}=[-1,1]$, and $j(0)=\min j = 0$. 
Here, $D(j)$ denotes the {\sl domain}\/
of the convex function $j$, i.e., the set 
where $j$ takes finite values.

It is well known that the graph $\beta$ 
induces maximal monotone operators (identically
noted as $\beta$ for simplicity) both in $H$ 
and in $L^2(Q)$, where $Q:=(0,T)\times\Omega$. 
For instance, one has
$\xi \in \beta(u)$ in the $H$-sense if and only if
$u,\xi \in H$ and $\xi(x) \in \beta(u(x))$
for a.e.~$x\in \Omega$. Moreover, let us define
the {\sl convex}\/ functional
\begin{equation}\label{defiJ}
  J:H\to [0,+\infty], \quad
   J(u):= \io j(u)\,\dix, 
\end{equation}
where the integral may well be $+\infty$ in the case when
$j(u)\not\in L^1(\Omega)$ (i.e., when $u\not\in D(J)$).
Then, $\beta = \de J$ in $H$, namely 
the operator induced by $\beta$ on $H$ coincides with
the $H$-subdifferential of the convex
functional $J$. As is customary when dealing with
multivalued operators, we shall often identify maximal
monotone operators with their graphs (cf., 
e.g., \cite{Ba,Br}).
With the above notation, equation~\eqref{eqn:intro},
where the coefficients $\epsi$ and $\delta$ have been set to 
$1$ for simplicity, becomes
\begin{equation}\label{eqn}
  u_{tt} 
   + A u_t
   + A u
   + \beta(u) 
   - \lambda u 
  \ni g
  \quext{in }\,(0,T)\times \Omega.
\end{equation}
Note the occurrence of the inclusion sign, motivated by
the fact that $\beta$ may be multi-valued.

In view of \eqref{Aneum} (or of \eqref{Adir}), \eqref{eqn}
can be read as a relation holding in $L^2(0,T;H)$ (and thus 
interpreted as a pointwise inclusion
almost everywhere in~$Q$). Indeed, \eqref{eqn}
looks as the most natural and appropriate weak formulation
of the strongly damped wave equation in the case when
$\beta$ is a smooth monotone function defined on the
whole real line. On the other hand, though \eqref{eqn}
is still perfectly meaningful from the mathematical viewpoint
under our assumptions \eqref{hp:f}-\eqref{hp:beta},
proving existence of solution in the current setting seems to be 
out of reach (see Remark~\ref{rem:reg1} below
for a counterexample in the spatially homogeneous case).
Mainly, what seems to fail is the possibility to interpret
point-by-point the equation, and in particular 
the constraint~$\beta$. 

Hence, we need to construct a furtherly relaxed formulation 
of the equation, for which one might be able to get existence. 
In performing this program, we would like
our new concept of solution to be still somehow
physically consistent. Namely, weak solutions
should comply with thermodynamical principles (like the energy inequality),  
satisfy a proper form of the constraint,
and be obtained as limit points of families of functions 
solving physically sound regularizations 
of the equation. To start with this program, we set
\begin{equation}\label{calV}
  \calV:= H^1(0,T;H)\cap L^2(0,T;V).
\end{equation}
Note that, in view of standard results on 
vector-valued functions, the above space coincides
with $H^1(Q)$ in the Neumann case. The duality
pairing between $\calV'$ and $\calV$ will be noted
by $\duavw{\cdot,\cdot}$. 
We also consider the space $\calH:=L^2(Q)=\LDH$ 
endowed with the natural scalar product, noted here as
$(\!( \cdot, \cdot )\!)$. Thanks to
standard results on Sobolev spaces, the inclusions 
$\calV\subset \calH \subset \calV'$ hold continuously 
and densely provided $\calH$ is identified with its 
dual by means of the above scalar product. Actually,
the weak formulation of our problem will 
strongly rely on the {\sl parabolic Hilbert triplet}~$(\calV,\calH,\calV')$.

We also need similar concepts in the case when
the time interval $(0,T)$ is replaced by 
$(0,t)$  for $0 <t\le T$. Namely,
%
%
we set $Q_t:= (0,t) \times \Omega$, and,
correspondingly, we note by $(\!( v, z )\!)\zt$ the 
(standard) scalar product in $\calH_t:=L^2(Q_t)$ and by
$\duavw{\cdot,\cdot}\zt$ the duality between
$\calV_t:= H^1(0,t;H)\cap L^2(0,t;V)$ and its dual.
We also set
\begin{equation}\label{Vtzero}
  \Vtzero:= \big\{ v\in \calV_t:~v\equiv 0~\text{on }\{t\}\times \Omega\big\},
\end{equation}
where relation $v\equiv 0$ is intended in the sense of traces (in time).
Clearly, $\Vtzero$ is a closed subspace of $\calV_t$.
Then, if $\fhi\in \Vtzero$,
extending it by $0$ for times larger than
$t$, we obtain an element of $\calV$, noted
in the following as $\tilfhi$. Correspondingly, if $\eta\in \calV'$,
we can naturally define its {\sl restriction}\/ $\eta^t$ 
to the time interval $(0,t)$ by setting, for $\fhi\in \Vtzero$,
\begin{equation}\label{restr}
  \duavw{ \eta^t,\fhi }\zt := 
   \duavw{ \eta,\tilfhi }.
\end{equation}
Actually, it is readily checked that $\eta^t\in \Vtzero'$. Moreover, the 
restriction operator $\eta\mapsto \eta^t$ is linear and continuous
from $\calV'$ to $\Vtzero'$.

With the above notation at disposal, we extend
the functional $J$ to time-dependent functions
by setting (see~\eqref{defiJ})
\begin{equation}\label{deficalJ}
  \calJ:\calH \to [0,+\infty], \quad
   \calJ(u):= \iTT J(u)\,\dit = \iTo j(u)\,\dix\,\dit,
\end{equation}
where, as before, the integral may also take
the value $+\infty$. 
Analogously, for $t\in(0,T]$, we put
\begin{equation}\label{deficalJt}
  \calJ_{(t)}:\calH_t \to [0,+\infty], \quad
   \calJ_{(t)}(u):= \itt J(u)\,\dis = \ito j(u)\,\dix\,\dis.
\end{equation}
As noted above, the $\calH$-subdifferential 
$\de \calJ(u)$ (or the analogue for $\calJ_{(t)}(u)$) 
can be still interpreted in the ``pointwise'' form $\beta(u)$.

We are now ready to introduce the weak form of the constraint $\beta$. 
We shall present most of the construction by working
on the time interval $(0,T)$. The adaptation to 
subintervals $(0,t)$ is straighforward and we mostly
leave it to the reader because we do not want
to overburden the notation. That said, we start by
setting $\calJ_{\calV}:=\calJ|_{\calV}$. 
It is readily proved that $\calJ_{\calV}$ is 
convex and lower semicontinuous on $\calV$.
Hence, we may take its subdifferential 
with respect to the duality pairing between 
$\calV$ and $\calV'$. Namely, for 
$\xi \in \calV'$ and $u\in \calV$, we put 
\begin{equation}\label{defibw}
  \xi \in \beta_w(u) \stackrel{\text{def}}{\Longleftrightarrow}
   u\in\calV,~~\xi\in \calV',~~\text{and }\,
   \duavw{\xi,v-u} + \calJ_{\calV}(u) \le \calJ_{\calV}(v) 
   ~~\text{for all }\, v \in \calV.
\end{equation}
%
%
%
The idea of ``relaxing'' $\beta$ in this way is not new;
for instance, the same method has been applied
in \cite{BCGG,DK,PS1} in other contexts. It is worth 
noting from the very beginning that $u\in D(\beta_w)$
still implies $u\in[-1,1]$ almost everywhere; in other words,
the weak operator $\beta_w$ still forces 
$u$ to assume only ``physically meaningful'' values.
Note that an alternative, but essentially 
equivalent, approach based on variational
inequalities has been devised in \cite{MZ} for the Cahn-Hilliard
equation with dynamic boundary conditions. The novelty
occurring in our case is related to the use
of ``parabolic'' (Sobolev-Bochner) spaces. Indeed,
this choice seems particularly appropriate for the 
present problem as far as it permits us to 
overcome some issues related with the (expected) 
low regularity of weak solutions.

Let us now characterize a bit more precisely the operator
$\beta_w$. We follow here
the lines of \cite{brezisart,GR} (see also \cite{BCGG}). Firstly, we observe
that (see, e.g., \cite[Prop.~2.3]{BCGG}), if 
$u\in \calV$, $\xi \in \calH$, and 
$\xi\in \beta(u)$ a.e.~in~$Q$, then $\xi\in \beta_w(u)$. 
Namely, if $\beta|_{\calV}$ denotes the 
restriction to $\calV$ of the ``pointwise'' operator~$\beta$,
then $\beta_w$ extends $\beta|_{\calV}$.
In other words, the ``strong'' constraint implies the ``weak'' one. 
Moreover (cf.~\cite[Prop.~2.5]{BCGG}), 
\begin{equation} \label{betaw2}
  \text{if }\,u\in \calV~\,\text{and }\,
  \xi \in \beta_w(u) \cap \calH, 
   \quext{then }\,\xi \in \beta(u)
   ~\,\text{a.e.~in }\,Q.
\end{equation}
In general, however, the elements $\xi\in\beta_w(u)$ 
(which lie, by definition, in the space $\calV'$) need
not belong to $\calH$. Hence,
the graph inclusion $\beta|_{\calV} \subset \beta_w$ 
is generally a proper one. 
Nevertheless, if $\xi \in \beta_w(u)$,
then $\xi$ ``automatically'' gains some more regularity.

In order to explain this phenomenon, we proceed along the lines of 
\cite[Sec.~2]{PS1}. Namely, for $t\in(0,T]$, we set 
\begin{subequations} \label{calX}
\begin{align} \label{calXneum}
  & \calX_t:=C^0(\ov{Q_t}),
  \quext{for Neumann boundary conditions},\\
 \label{calXdir}
  & \calX_t:=\big\{ u\in C^0(\ov{Q_t}):~u\equiv 0~
   \text{on }\,[0,t]\times \Gamma\big\},
 \quext{for Dirichlet boundary conditions.}
\end{align}
\end{subequations}
For $t=T$ we simply write $\calX=\calX_T$.
We also set, in both cases,
\begin{equation} \label{calX0}
  \Xtzero :=
   \big\{ v\in \calX_t:~v\equiv 0~\text{on }\{t\}\times \Omega\big\}.
\end{equation}
The space $\calX_t$ (hence its closed subspace $\Xtzero$)
is naturally endowed with the supremum norm $\| \cdot \|_\infty$.
Moreover, also thanks to
the smoothness of $\Omega$ in the Neumann
case, $\calX_t\cap \calV_t$ is dense both in $\calX_t$ 
and in $\calV_t$. Let now $\xi\in \calV'$ (the analogue applies 
with straighforward modifications to $\xi\in \calV_t'$) and let us suppose
that $\xi$, if restricted to the functions $\fhi\in\calX\cap \calV$,
is continuous with respect to the $\calX$-norm, i.e., 
there exists $C>0$ such that 
\begin{equation} \label{repres}
  \big| \duavw{\xi,z} \big| \le C \| z \|_{\infty}
   \quext{for any~\,$z\in \calX \cap \calV$}.
\end{equation}
In that case, by density, $\xi$ extends {\sl in a unique way}\/
to a bounded linear functional on $\calX$. Namely, there exists a unique
$\calT\in \calX'$, which can be seen as a Borel measure on $\ov{Q}$
in view of Riesz' representation theorem, such that
\begin{equation} \label{identif}
  \duavw{\xi,z}=\iint_{\ov{Q}} z\,\diT
   \qquext{for any~\,$z\in \calX \cap \calV$}.
\end{equation}
In this situation we say that the measure $\calT$ {\sl represents}~$\xi$ 
on~$\calX$. Actually this situation automatically
occurs when $\xi$ is an element of a weak constraint. Indeed,
by an easy adaptation of~\cite[Prop.~2.1]{PS1}
(which, in turn, is based on the results of~\cite{brezisart}),
one can see that, up to some 
adjustment related to the boundary
conditions, any $\xi\in\beta_w(u)$, when 
restricted to continuous functions, is 
represented by a measure $\calT$ 
defined on the parabolic cylinder $Q$.
Such a measure, in turn, is related to the original
operator $\beta$ in the following way (cf.~\cite[Thm.~3]{brezisart} 
for further details): noting as  
$\calT=\calT_a+\calT_s$ the Radon-Nikodym 
decomposition of $\calT$,
with $\calT_a$ ($\calT_s$, respectively) 
standing for the absolutely continuous 
(singular, respectively) part, we then have
\begin{align}\label{bre1}
  & \calT_a u\in L^1(Q),\\
 \label{bre2}
  & \calT_a(t,x) \in \beta(u(t,x)) \quext{for a.e.~$(t,x)\in Q$,}\\
 \label{bre3}
  & \duavw{\xi,u} - \iTo \calT_a u\,\dix\,\dit
    = \sup \bigg\{\iint_{\ov{Q}} z\,\deriv\!\calT_s,~z\in \calX,~
    z(\ov{Q}) \subset [-1,1] \bigg\}.
\end{align}
Hence, the continuous part $\calT_a$ of the measure~$\calT$
satisfies the constraint pointwise
(in view of \eqref{bre2}), whereas the singular part 
$\calT_s$ is characterized by \eqref{bre3}. 

In particular, we expect that condition \eqref{bre3} 
could be made more precise. Namely,
noting as $\calT_s = \rho |\calT_s|$ the {\sl polar 
decomposition}\/ of $\calT_s$, where $|\calT_s|$ is 
the total variation of $\calT_s$,
proceeding along the lines of \cite[Thm.~3]{GR}
one may prove that
\begin{equation} \label{bre4}
  \rho \in \de I_{[-1,1]}(u)
   ~~\text{$|\calT_s|$-a.e.~in }\,\ov{Q}.
\end{equation}
In other words, we expect the singular part of $\calT$ to
be supported on the set where $|u|=1$ and that $\rho=1$ 
where $u=1$, $\rho=-1$ where $u=-1$. In this sense, 
also the singular part of $\calT$ is, at least partially,
reminiscent of the expression of the graph $\beta$.

Actually, the characterization \eqref{bre4} is proved in \cite{GR}
in the case when $\calV=H^1_0(\Omega)$, $\Omega$ a bounded domain of
$\RR^N$, and may be likely extended 
to the present situation. However,
a detailed proof  may involve some technicalities 
particularly related to the facts that we are working in 
the parabolic cylinder and should distinguish between
the Dirichlet and Neumann cases. For this reason,
we omit details here. 
%
%
We note, however, that \eqref{bre4} is straighforward whenever
we additionally know that $u\in \calX\cap \calV$ (i.e., $u$,
beyond lying in $\calV$, is continuous).
Indeed, in that case, from \eqref{bre3} there follows
\begin{align} \label{bre3b}
  \iint_{\ov{Q}} \rho u \,\deriv\!|\calT_s|
   & = \iint_{\ov{Q}} u \,\deriv\!\calT_s
     = \duavw{\xi,u} - \iTo \calT_a u \,\dix\,\dit \\
 \no
  & = \sup \bigg\{\iint_{\ov{Q}} z\,\deriv\!\calT_s,~z\in \calX,~
    z(\ov{Q}) \subset [-1,1] \bigg\}
    = | \calT_s |(\ov{Q}),  
\end{align} 
the latter term denoting the total variation of the measure
$\calT_s$. Comparing terms, we then deduce
$\rho u=1$ $|\calT_s|$-a.e.~in~$\ov{Q}$, as desired.
\beos\label{rem:traccia}
 It is worth observing that, in the Neumann
 case, the singular component $\calT_s$
 of the measure $\calT$ representing $\xi\in\beta_w(u)$ may 
 be, at least partially, supported on the boundary of~$Q$.
 Let us see this by a simple one-dimensional example. 
 Let $Q=(-1,1)$, $u(x)=x$, $j=I_{[-1,1]}$. 
 Then, if $\xi=\alpha\delta_1$ for some $\alpha>0$,
 where $\delta_1$ is the Dirac delta concentrated in~$1$,
 it is clear that, for any $v\in H^1(-1,1)$ such that 
 $v\in D(J)$ (i.e., such that $-1\le v(x)\le 1$ for 
 all $x\in Q$), there holds
 \begin{equation} \label{bre3c}
   \duav{\xi,v-u} 
   = \alpha (v(1) - u(1))
   = \alpha (v(1) - 1)
   \le 0 
   = J(v) - J(u).
 \end{equation} 
 Hence, $\xi\in \beta_w(u)$ by definition of subdifferential.
\eddos
\noindent%
Having clarified the nature of the weak constraint $\beta_w$, we 
can now observe that equation~\eqref{eqn} admits a natural
energy functional
\begin{equation}\label{defiE}
  \calE(u,u_t) 
   = \io \Big(
    \frac12 | u_t |^2 
    + \frac12 | \nabla u |^2
    + j(u) - \frac\lambda2 u^2 \Big)\,\dix
  = \frac12 \big( \| u_t \|^2 
    + \| \nabla u \|^2
    + 2 J(u) - \lambda \| u \|^2 \big).
\end{equation}
Indeed, testing \eqref{eqn} by $u_t$ and integrating 
in time and space, one can get that the value of~${\cal E}$ 
at any time $t>0$ is bounded by the initial value~${\cal E}(u_0,u_1)$ 
plus the power of the external applied forces $g$ (see 
\eqref{enerlim2} below).   
Actually, as will be explained later on, the low regularity 
of solutions does not allow us to perform this estimate
directly for weak solutions, but only for a suitable approximation
of the problem. This is basically the reason for which we 
will only be able to prove an energy {\sl inequality}\/
for weak solutions, cf.~Theorem~\ref{teo:esi} below.

We can now introduce our assumption on the 
source term and on the initial data, the latter
corresponding exactly to the finiteness of the 
``initial energy'':
\begin{align} \label{hp:g}
  & g \in \LDH, \\
  \label{hp:u0}
  & u_0 \in V, \quad u_1\in H, \quad J(u_0)<\infty.
\end{align}
Then, we can make precise our concept of weak 
solution (to be precise, we shall speak of ``parabolic
duality weak solution'' or something similar, but we will
rather use ``weak solution`` just for simplicity):
\bede\label{def:sol}
 A couple $(u,\xi)$ is called a~{\rm weak solution} to
 the initial-boundary value problem for 
 the strongly damped wave equation with constraint 
 whenever the following conditions hold:
 \begin{itemize}
 \item[(a)] There hold the regularity properties
 \begin{align}\label{regou}
   & u \in W^{1,\infty}(0,T;H) \cap  H^{1}(0,T;V),\\
  \label{regoxi}
   & \xi \in \calV'\cap\calX';
 \end{align}
 more precisely, by \eqref{regoxi} we may require 
 that there exists a measure $\calT\in \calX'$ representing 
 $\xi$ over $\calX$ in the sense of~\eqref{identif}. 
 \item[(b1)] There holds the following weak version 
 of~\eqref{eqn}:
 \begin{align}\no
   & \mbox{} - (\!( u_t, \fhi_t )\!) 
    + ( u_t(T) , \fhi(T) )
    + (\!( \nabla u_t, \nabla \fhi )\!) 
    + (\!( \nabla u, \nabla \fhi )\!)\\
  \label{eqnwT}  
   & \mbox{} ~~~~~
    + \duavw{ \xi,\fhi }
    - \lambda (\!( u, \fhi )\!) 
    = ( u_1 , \fhi(0) )
    + (\!( g, \fhi )\!) , \quad  \forall\,\fhi\in\calV.
 \end{align}
 \item[(b2)] An analogue of~\eqref{eqnwT} holds also on subintervals,
 in the following sense: for any $t\in (0,T]$ there exists a 
 functional $\xi_{(t)}\in \calV_t'$ such that 
 \begin{align}\no
   & \mbox{} - (\!( u_t, \fhi_t )\!)\zt 
    + ( u_t(t) , \fhi(t) )
    + (\!( \nabla u_t, \nabla \fhi )\!)\zt 
    + (\!( \nabla u, \nabla \fhi )\!)\zt\\
  \label{eqnw}  
   & \mbox{} ~~~~~
    + \duavw{ \xi_{(t)},\fhi }\zt- \lambda (\!( u, \fhi )\!)\zt 
    = ( u_1 , \fhi(0) )
    + (\!( g, \fhi )\!)\zt ,\quad  \forall\,\fhi\in\calV_t.
 \end{align}
 Moreover, $\xi_{(t)}$ lies in $\calV_t'\cap\calX_t'$; hence
 it is (uniquely) represented over $\calX_t$
 by a measure $\calT_{(t)}$. In addition to that, for every 
 $t\in (0,T]$, the functionals $\xi$ and 
 $\xi_{(t)}$ are\/ {\rm compatible}, namely,
 for every $\fhi\in\Vtzero$ we have 
 ($\tilfhi$ denoting the trivial extension of $\fhi$)
 \begin{equation}\label{incl}
   \duavw{ \xi_{(t)},\fhi }\zt 
    = \duavw{ \xi,\tilfhi }. 
 \end{equation}
 In other words, the functional $\xi_{(t)}$, when computed 
 on the elements of $\Vtzero$, coincides with the canonical
 restriction $\xi^t$ of $\xi$ (cf.~\eqref{restr}).
 \item[(c)] There holds the inclusion
 \begin{equation}\label{incl2}
    \xi \in \beta_w(u) 
    \quext{in }\,\calV'.
 \end{equation}
 More generally, for every $t\in (0,T]$, $\xi_{(t)}\in \beta_{w,(t)}(u)$
 in $\calV_t'$. Here $\beta_{w,(t)}$ represents the weak version of $\beta$ 
 in the interval $(0,t)$; namely, $\beta_{w,(t)}$ is the subdifferential 
 of the restriction of $\calJ_{(t)}$ to $\calV_t$ with respect to the duality
 product between $\calV_t$ and $\calV_t'$.
 \item[(d)] There holds the Cauchy condition 
 \begin{equation}\label{init}
   u|_{t=0}=u_0.
 \end{equation}
 \item [(e)] For every $s,t\in [0,T]$, the couple $(u,\xi)$ satisfies the equality
 \begin{align}\no
   & \mbox{} 
    - \| u_t \|_{L^2(s,t;H)}^2
    + ( u_t(t) , u(t) )
    + \frac12 \| \nabla u(t) \|^2
    + \| \nabla u \|_{L^2(s,t;H)}^2 
    + \duavw{ \xi_{(t)}, u }_{(0,t)}
    - \duavw{ \xi_{(s)}, u }_{(0,s)} \\
  \label{ener}
   & \mbox{} ~~~~~
    - \lambda \| u \|_{L^2(s,t;H)}^2
   = ( u_t(s) , u(s) )
    + \frac12 \| \nabla u(s) \|^2
    + \int_s^t \io g u\,\dix\,\ditau.
 \end{align}
 \end{itemize}
\edde
\noindent
It is worth discussing a bit how the above formulation has been obtained
from~\eqref{eqn}. First of all, $\beta$ has been replaced with its ``relaxed'' 
form $\beta_w$. Correspondingly, \eqref{eqn} has been
restated in the ``parabolic'' dual space $\calV'$ by using the
test function $\fhi\in \calV$ and performing suitable 
integrations by parts. In particular, a key 
point stands in the integration in time of the ``hyperbolic''
term $u_{tt}$. Indeed, no second time derivatives 
of $u$ appear in \eqref{eqnwT} (or in \eqref{eqnw}). 
In addition to that, the 
Cauchy condition for $u_t$ is now ``embedded'' 
into~\eqref{eqnwT} and \eqref{eqnw}.
\beos\label{rem:tempo}
 We need to explain in some detail the ``meaning'' of~\eqref{eqnw},
 especially in relation with the constraint term. Actually, if $\xi\in \calV'$,
 there is no canonical way of restricting $\xi$ 
 to obtain an element of $\calV_t'$. The best we can do is
 restricting $\xi$ as explained in \eqref{restr} to obtain
 a functional $\xi^t\in\Vtzero'$. However, writing 
 \eqref{eqnw} as a relation in $\Vtzero'$ (i.e., considering
 only test functions $\fhi\in \Vtzero$) would give rise 
 to some information loss. Namely, it may happen that
 the singular part $\calT_s$ of $\calT$ 
 is, at least partially, supported on 
 some set of the form $\{t\}\times\ov{\Omega}$
 (or, correspondingly, $\duavw{\xi_{(t)},\fhi}_{(0,t)}$  
 may also depend on the trace of $\fhi$
 on $\{t\}\times\ov{\Omega}$).
\eddos
\beos\label{rem:reg1}
 It is worth noting that, according to the above definition,
 $u_t$ need not be continuous with respect to time,
 independently of the target topology. This fact is 
 a distinctive feature of this problem and there seems to be
 no hope of avoiding jumps of $u_t$, at least for a {\sl general}\/
 constraint $\beta$. Here is a simple example where 
 a jump occurs. Let us consider the case of
 spatially homogeneous solutions to the Neumann problem
 with $\lambda = 0$ and $g = 0$. In other words, we reduce our 
 problem to the ``toy model''
 represented by the ODE
 \begin{equation}\label{toy}
   u_{tt} + \beta(u) \ni 0,
 \end{equation}
 a weak solution $u$ to which exists according to our theory. 
 Let us also choose $\beta=\de I_{[-1,1]}$. Then,  
 if we take, for instance, $u_0=0$ and $u_1=1$, we get
 that $u(t)=t$ at least for $t\in[0,1)$.
 As $t$ gets to $1$, $u_t$ {\sl must}\/
 develop a discontinuity, otherwise,
 $u(t)$ would become {\sl strictly}\/ larger than~$1$
 for $t>1$, and the equation would no longer
 make sense. Hence, the only possibility for 
 the trajectory $u(t)$ is to jump instantaneously in such a way
 that, in a right neighbourhood $(1,1+\epsilon)$ of $t=1$,
 \begin{equation}\label{counter}
   u_t(t) = \ell, \quad
    u(t) = 1 + \ell ( t - 1 ), \quad \ell \le 0.
 \end{equation}
 The trajectory, at least in principle, 
 may ``choose'' at which level $\ell$ the time
 derivative $u_t$ ``decides'' to jump (hence we have no uniqueness). 
 If it jumps to $\ell<0$, then $u(t)$ starts to decrease 
 from the value $1$ at a constant velocity $\ell$
 until it reaches the value $-1$ (where a new jump of 
 $u_t$ must occur). On the contrary, if
 $u_t$ jumps to $\ell=0$, then it will be either
 $u(t)=1$ and $u_t(t)=0$ forever, or 
 after some time $u_t$ may make a further jump to
 some $\ell<0$, starting from which $u$ begins 
 to decrease as specified above. More precisely, we can
 notice that, for \eqref{toy}, the weak formulation
 over~$(0,T)$ (cf.~\eqref{eqnwT}) reads
 \begin{equation}\label{toyweak}
   - \int_0^T u_{t} \fhi_t \,\dit 
    + \fhi(T) u_t(T)
    + \duavw{ \xi, \fhi }
    = \fhi(0) u_1
   \quad \perogni \fhi\in H^1(0,T).
 \end{equation}
 Hence, it is easy to check that, for all $\ell\le 0$,
 the function $u$ described above solves~\eqref{toyweak}
 on a suitable interval $[0,T]$ with $T>1$ chosen sufficiently
 small so that no other jumps of $u_t$ occur.
 Note in particular that different 
 choices of $\ell$ correspond to different ``values''
 of $\xi\in\beta_w(u)$. Indeed, from \eqref{toyweak}
 we get
 \begin{equation}\label{toyweak2}
   - \int_0^1 \fhi_t \,\dit 
   - \int_1^T \ell \fhi_t \,\dit 
    + \ell \fhi(T)
    + \duavw{ \xi, \fhi }
    = \fhi(0),
 \end{equation}
 whence 
 \begin{equation}\label{toyweak3}
   \duavw{ \xi, \fhi } = (1 - \ell) \fhi(1), 
 \end{equation}
 or, in other words, $\xi=(1-\ell)\delta_{t=1}$
 ($\delta$ standing for the Dirac delta)
 and we can notice that this is consistent with the above 
 characterization of $\beta_w$. Indeed, at least for 
 $\ell<0$, $t=1$ is the only time at which 
 $u$ takes the value $1$ and $\xi$ may have, 
 and in fact has, a ``singular''
 part. However, we will see in the sequel 
 (cf.~Remarks~\ref{energia2} and \ref{energia3} below) 
 that not every jump of $u_t$ (or, in the current 
 example, every value of $\ell>0$) is ``physically''
 admissible. 
 %
 %
\eddos
\noindent%
We can now introduce the statement of 
our main result:
\bete\label{teo:esi}
 Let us assume \eqref{hp:f}, \eqref{hp:beta}, \eqref{hp:g}, 
 and \eqref{hp:u0}. Then, there exists at least one weak solution 
 $(u,\xi)$, in the sense 
 of\/~{\rm Definition~\ref{def:sol}}, to
 the initial-boundary value problem for 
 the strongly damped wave equation with constraint. 
 Moreover $u_t\in BV(0,T;X)$ for any Banach space $X$ 
 such that $L^1(\Omega)$ and $V'$ are compactly embedded in $X$.\\[1mm]
 In addition, for almost every $s\in [0,T)$ (surely including
 $s=0$) and every $t\in (s,T]$,
 the following version of the energy inequality holds:
 \begin{equation}\label{enerlim2}
  \calE(u(t),u_t(t)) 
   + \| \nabla u_t \|_{L^2(s,t;H)}^2
  \le \calE(u(s),u_t(s)) 
   + \int_s^t ( g,u_t )\,\ditau,
 \end{equation}
 where $\calE$ is defined in~\eqref{defiE}.\\[1mm]
 Finally, in the case when we additionally have  
 \begin{equation}\label{hp:u0:2}
   u_0 \in D(A),
 \end{equation}
 then $u$ enjoys the additional regularity property
 \begin{equation}\label{regou2}
   u \in C_w([0,T];D(A)).
 \end{equation}
 Namely, $u(t)$ belongs to $D(A)$ for every $t\in[0,T]$ and 
 $t\mapsto u(t)$ is continuous when the target space is endowed
 with the weak topology.
\ente
%


\section{Proof of Theorem~\ref{teo:esi}}
\label{sec:proofs}

\textbf{Step 1. Approximation.}~~%
We start by introducing a natural regularization of (the strong form
of) equation~\eqref{eqn} depending on an approximation parameter $\epsi$ 
(which will then be let go to $0$). To this aim, for $\epsi\in(0,1)$, we let 
$j\ee:\RR\to[0,\infty)$ denote the {\sl Moreau-Yosida regularization}\/
of $j$ (cf., e.g., \cite{Br} for details). In particular, $j\ee$ turns out to be   
convex and lower semicontinuous. Moreover, its derivative
$\beta\ee:=\de j\ee$ corresponds to the {\sl Yosida approximation}\/
of $\beta=\de j$. Under our assumptions $\beta\ee$ is monotone and globally 
Lipschitz continuous on the whole real line and it satisfies
$\beta\ee(0) = 0$. We also set
\begin{equation}\label{defiJee}
  J\ee(u):= \io j\ee(u)\,\dix, \quad 
   \calJ\ee(u):= \iTo j\ee(u)\,\dix\,\dit.
\end{equation}
Moreover, we regularize the initial data by
taking, for $\epsi\in(0,1)$, 
$u\ee_0$ and $u\ee_1$ satisfying
\begin{align}\label{inizee1}
  & u\ee_0 \in D(A), \quad u\ee_1 \in V,
   \quad J\ee(u\ee_0) \le J(u_0), \\
 \label{inizee2}
  & u\ee_0 \to u_0~~\text{in }\, V, \qquad
   u\ee_1 \to u_1~~\text{in }\, H.
\end{align}
The construction of approximate initial data complying with
\eqref{inizee1}-\eqref{inizee2} is standard. For instance, 
one may take $u\ee_0$ as the solution to the elliptic 
singular perturbation problem
\begin{equation}\label{ellee}
  u\ee_0\in D(A), \qquad
   u\ee_0 + \epsi A u\ee_0 = u_0~~\text{in }\,H.
\end{equation}
In particular, the last of~\eqref{inizee1}
can be shown by testing the equation
in~\eqref{ellee} by $\beta\ee(u\ee_0)$
and noting that
\begin{equation}\label{ellee2}
  ( \beta\ee(u\ee_0), u\ee_0 - u_0 )
   \ge J\ee(u\ee_0) - J\ee(u_0)
   \ge J\ee(u\ee_0) - J(u_0),
\end{equation}
the latter inequality following from the monotonicity
of the Moreau-Yosida regularization $J\ee$ with respect to $\epsi$.

We are now ready to introduce our approximated equation
\begin{equation}\label{eqne}
  u_{tt}\ee 
   + A u_{t}\ee
   + A u\ee
   + \beta\ee(u\ee) 
   - \lambda u\ee 
  = g
  \quext{in }\,(0,T)\times \Omega.
\end{equation}
Correspondingly, we have the following well-posedness
and regularity result:
\bete\label{exi:e}
 Let us assume \eqref{hp:f}, \eqref{hp:beta}, \eqref{hp:g}, 
 and \eqref{hp:u0}. For $\epsi\in(0,1)$, let
 $u_0\ee$, $u_1\ee$, and $\beta\ee$ be as detailed
 above. Then, there exists a unique solution 
 \begin{equation}\label{regouee}
   u\ee \in H^2(0,T;H) \cap W^{1,\infty}(0,T;V) \cap H^1(0,T;D(A))
 \end{equation} 
 to equation \eqref{eqne}, complemented with the initial conditions 
 $u^\epsi|_{t=0}=u_0\ee$ and $u_t^\epsi|_{t=0}=u_1\ee$.
 Moreover, for all $t,s\in[0,T]$, the following energy\/
 {\rm equality} holds:
 \begin{align}\no
   & \mbox{} 
    \frac{1}{2}\| u^\epsi_t(t) \|^2
    + \int_\Omega j\ee(u^\epsi(t)) \,\dix
    - \frac{\lambda}{2} \|u^\epsi(t)\|^2
    + \frac{1}{2}\| \nabla u^\epsi (t)\|^2
    + \| \nabla u^\epsi_t \|_{L^2(s,t;H)}^2\\
  \label{enereps}
   & \mbox{} ~~~~~
    =\frac{1}{2}\| u^\epsi_t(s) \|^2
    + \int_\Omega j\ee(u^\epsi(s))\,\dix
    - \frac{\lambda}{2} \|u^\epsi(s)\|^2
    + \frac{1}{2}\| \nabla u^\epsi (s)\|^2
    + \int_s^t ( g,u^\epsi_t )\,\ditau.
\end{align}
\ente
\noindent%
The proof of Theorem~\ref{exi:e} is fairly standard and could be carried out, e.g.,
by following the lines of~\cite{LLM}. Here it is just worth noting that 
the regularity conditions stated in~\eqref{regouee} are compatible with 
the assumptions \eqref{inizee1}-\eqref{inizee2} on the regularized initial
data. Moreover, one could easily check that~\eqref{regouee} can be 
(at least formally) obtained testing \eqref{eqne} 
by $u\ee_{tt} + A u\ee_t$ and performing
integrations by parts. In particular, the term $\beta\ee(u\ee)$
can be managed thanks to the Lipschitz continuity of Yosida approximations.
It is also worth noting that, in this regularity setting, equation 
\eqref{eqne} makes sense pointwise; indeed, all its single terms belong
to the space $\LDH$ (in particular, we do not need to regularize 
the source term $g$: condition \eqref{hp:g} is enough).
Hence, testing the equation by $u\ee_t$ is allowed:
this gives relation \eqref{enereps} by means of well-known 
chain rule formulas.

\medskip

\noindent%
\textbf{Step 2. A priori estimates.}~~%
We now derive a number of bounds, uniform with respect to 
the regularization parameter~$\epsi$,
for the solutions $u\ee$ given by Theorem~\ref{exi:e}.
First of all, testing \eqref{eqne}
by $\fhi\in\mathcal V_t$, 
integrating over $Q_t$, and performing suitable integrations
by parts (both in space and in time), we deduce the
integrated (weak) formulation
\begin{align}\no
   & \mbox{} - (\!( u^\epsi_t, \fhi_t )\!)\zt 
    + ( u^\epsi_t(t) , \fhi(t) )
    + (\!( \nabla u^\epsi_t, \nabla \fhi )\!)\zt
    + (\!( \nabla u^\epsi, \nabla \fhi )\!)\zt\\
  \label{eqnweps}  
   & \mbox{} ~~~~~
    + (\!( \beta\ee(u^\epsi),\fhi )\!)\zt - \lambda(\!( u^\epsi, \fhi )\!)\zt
    = ( u_1\ee , \fhi(0) )
    + (\!( g, \fhi )\!)\zt,
   \quad \perogni \fhi\in\calV_t.
\end{align}
Of course, \eqref{eqnweps} holds in particular for $t=T$ and $\fhi\in \calV$.
Next, setting $s=0$ in \eqref{enereps}, or, in other words, testing \eqref{eqne}
by $u_t\ee$ and integrating over~$Q_t$, $t\in (0,T]$, we find
(by Young's inequality) 
\begin{align}\no
 & \mbox{} 
   \frac{1}{2}\| u^\epsi_t(t) \|^2
    + \int_\Omega j\ee(u^\epsi(t)) \,\dix
    - \frac{\lambda}{2} \|u^\epsi(t)\|^2
    + \frac{1}{2}\| \nabla u^\epsi (t)\|^2
    +  \| \nabla u^\epsi_t \|_{L^2(0,t;H)}^2
   = M_1(\epsi) + (\!(  g,u^\epsi_t )\!)_{(0,t)}\\
     & \mbox{} ~~~~~\label{Gronw}
     \leq M_1(\epsi) + \|g\|^2_{L^2(0,t;H)}
      + \frac{1}{4} \|u^\epsi_t\|^2_{L^2(0,t;H)}
   =: M_2(\epsi) + \frac{1}{4} \|u^\epsi_t\|^2_{L^2(0,t;H)},
\end{align}
where we have set
\begin{align}\label{rs11}
  & M_1(\epsi):=\frac{1}{2}\| u_1\ee \|^2
     + J\ee(u_0\ee)
     - \frac{\lambda}{2} \|u\ee_0\|^2
     + \frac{1}{2}\| \nabla u\ee_0\|^2,
       \quad
    M_2(\epsi):=M_1(\epsi) +\|g\|^2_{L^2(0,t;H)},
\end{align}    
and we may notice that actually $M_1$ and $M_2$ are bounded uniformly
in~$\epsi$ due to \eqref{inizee1}-\eqref{inizee2},
\eqref{hp:u0}, and \eqref{hp:g}.
Let us also observe that,
thanks to the properties of the Yosida approximation
(cf.~\cite{Br}), there exists a constant $c\ge 0$ such that 
$j\ee(r)- \lambda r^2\geq-c$ 
for all $r\in\RR$ and all $\epsi\in (0,1)$.
Hence, applying Gronwall's lemma to \eqref{Gronw},
we obtain
\begin{subequations}\label{est:eps}
 \begin{align}\label{est:eps:1}
     & \| u^\epsi_t(t) \| 
      \leq M\quad\mbox{for all }\,t\in[0,T],
 \end{align}
 and for all $\epsi\in(0,1)$. Here and below $M$ denotes a positive constant, 
 possibily different from line to line, depending on the problem data, 
 but independent of $\epsi$. From \eqref{Gronw} we also get
 \begin{align}
     &\| u^\epsi \|_{H^1(0,T;V)}\leq M, \label{est:eps:2}\\
     &\int_\Omega j\ee(u^\epsi(t))\,\dix \leq M\quad\mbox{for all }\,t\in[0,T],\label{est:eps:3}
 \end{align} 
 for all $\epsi\in (0,1)$.
\end{subequations}
Setting now $\fhi=u\ee$ in \eqref{eqnweps} and taking $t=T$, we deduce
\begin{align}\no
   & \mbox{} 
    \iTo \beta\ee(u^\epsi) u^\epsi\,\dix\,\dit 
     + \| \nabla u^\epsi \|_{L^2(0,T;H)}^2
     + \frac12 \| \nabla u^\epsi(T) \|^2
    = \| u^\epsi_t \|_{L^2(0,T;H)}^2
    - ( u^\epsi_t(T) , u^\epsi(T) ) \\
  \label{enereps2}
   & \mbox{} ~~~~~
   + \frac12 \| \nabla u_0^\epsi \|^2
   + \lambda \| u^\epsi \|_{L^2(0,T;H)}^2
   + ( u_1\ee , u_0\ee )
   + (\!( g, u^\epsi )\!).
\end{align}
Now, the \rhs\ is bounded uniformly~in~$\epsi$ due 
to \eqref{est:eps}, \eqref{hp:g}-\eqref{hp:u0} 
and \eqref{inizee1}-\eqref{inizee2}. Moreover, it 
is easy to check (cf.~also \cite[Appendix]{MZ})
that there exist constants $c_1>0$, $c_2\ge 0$ 
independent of $\epsi\in(0,1)$ such that 
$c_1 | \beta\ee (r)| \leq \beta\ee (r) r + c_2$ 
for all $r\in\RR$. Hence, \eqref{enereps2} entails
\begin{subequations}
\begin{align}\label{est:eps:4}
  & \|\beta\ee(u^\epsi)\|_{L^1(0,T;L^1(\Omega))}\leq M,
\end{align}
for all $\epsi\in(0,1)$. Then, using once more \eqref{est:eps}
and comparing terms in \eqref{eqne}, we also find
\begin{align}\label{est:eps:5}
  & \| u_{tt}^\epsi\|_{L^1(0,T;X)}\leq M,
\end{align}
for any Banach space $X$ such that $L^1(\Omega)\subset X$ 
and and $V'\subset X$ with continuous and compact embeddings.
In particular, it is not restrictive to assume $X$ be the 
dual of a reflexive and separable space
(for instance, in dimension $N=3$, one may take 
$X=H^{-2}(\Omega)$). Hence, we have 
\begin{align}
  & \|u_{t}^\epsi\|_{W^{1,1}(0,T;X)}\leq M, \label{est:eps:6}
\end{align}
\end{subequations}
for all $\epsi\in (0,1)$.

\medskip

\noindent%
\textbf{Step 3. Passage to the limit.}~~%
Now we aim at letting~$\epsi\searrow0$. From \eqref{est:eps:1}-\eqref{est:eps:2}
we deduce that there exists a function 
$u$ of the regularity specified in 
\eqref{regou} such that
\begin{subequations}\label{limit:eps}
\begin{align}
  & u^\epsi\rightharpoonup u \quad\mbox{weakly in }\,H^1(0,T;V),\label{limit:eps:1}\\
  & u^\epsi\rightharpoonup u \quad\mbox{weakly star in }\,W^{1,\infty}(0,T;H),\label{limit:eps:1*}
\end{align}
and in particular 
\begin{align}
  & u^\epsi(t)\rightarrow u(t) \quext{strongly in }\,H\, \mbox{ for all }t\in[0,T],\label{limit:eps:2a}\\
  & u^\epsi(t)\rightharpoonup u(t) \quext{weakly in }\,V\, \mbox{ for all }t\in[0,T]\label{limit:eps:2b}.
\end{align}
It is worth stressing that the above convergence relations, as well as the
ones that will follow, are intended to hold up to extraction of suitable (nonrelabelled)
subsequences of $\epsi\searrow0$. 

Since $V$ is compactly embedded into $H$, in view of conditions \eqref{est:eps:2} 
and \eqref{est:eps:6} we can apply 
\cite[Corollary 4]{Si} with the three spaces $V\subset\subset H\subset X$ and $p=2$ in order to obtain that
\begin{align}\label{limit:eps:4}
  u^\epsi_t\to u_t \quext{strongly in }\,L^2(0,T;H). 
\end{align}
Moreover, condition \eqref{est:eps:6} implies 
that the functions $u_t^\epsi$ are uniformly bounded in $BV(0,T;X)$
(for the properties of vector-valued $BV$-spaces one can 
refer, e.g., to \cite[Appendix]{Br}).
In view of the fact that we may assume $X$ be the dual of a 
reflexive and separable space,
we can employ a generalization of Helly's theorem 
(cf., e.g., \cite[Thm.~3.1]{MM}
or \cite[Lemma 7.2]{DMDSM}), 
providing a function $v\in BV(0,T;X)$ such that 
\begin{align}\label{limit:eps:2c}
  u_t^\epsi(t)\rightharpoonup v(t) \quext{weakly star in }X\, \mbox{ for all }t\in[0,T].
\end{align}
It is easily seen that $v$ coincides with $u_t$ almost everywhere.
Hence, up to changing the representative of $u_t$, we may assume $v=u_t$ everywhere 
on $[0,T]$. Moreover, combining \eqref{limit:eps:2c} with 
\eqref{est:eps:1}, we obtain
\begin{align}\label{limit:eps:3}
 u_t^\epsi(t)\rightharpoonup u_t(t) \quext{weakly in }H\,\mbox{ for all }t\in[0,T].
\end{align}
Moreover we get
\begin{align}\label{limit:eps:2c*}
  u_t^\epsi\rightharpoonup u_t \quext{weakly star in }\,BV(0,T;X).
\end{align}
Let us now show that the functions $\beta\ee(u^\epsi)$ are uniformly bounded
(with respect to $\epsi$) in~$\calV'$. Actually, writing \eqref{eqnweps}
for $t=T$, and using Holder's inequality, \eqref{inizee2}, \eqref{hp:g} and the estimates
\eqref{est:eps:1} and \eqref{est:eps:2}, we find  
\begin{align}\no
  & \big| \duavw{ \beta\ee(u^\epsi),\fhi } \big| 
    \leq \| u^\epsi_t\|_{L^2(0,T;H)} \|\fhi_t\|_{L^2(0,T;H)}
    + \| u^\epsi_t(T)\| \|\fhi(T)\| \\ 
 \no
  & \mbox{} ~~~~~~~~~~
   + \| \nabla u^\epsi_t\|_{L^2(0,T;H)} \|\nabla \fhi\|_{L^2(0,T;H)}
   + \| \nabla u^\epsi\|_{L^2(0,T;H)} \| \nabla \fhi\|_{L^2(0,T;H)} \\
  & \mbox{} ~~~~~~~~~~
   + \lambda\| u^\epsi\|_{L^2(0,T;H)} \|\fhi \|_{L^2(0,T;H)}
   + \| u_1^\epsi\| \| \fhi(0) \|
   + \| g \|_{L^2(0,T;H)} \|\fhi\|_{L^2(0,T;H)}
   \leq C \|\fhi\|_{\mathcal V}, \label{eta:inV}
\end{align}
for all $\fhi\in\mathcal V$, where $C>0$ is independent of $\epsi$. 
Therefore, we can infer that there exists $\xi\in \calV'$ such that
\begin{align}\label{limit:eps:5}
  & \beta\ee(u^\epsi)\rightharpoonup \xi\quext{weakly in }\,\mathcal V'.
\end{align}
Next, recalling the definition \eqref{calX} of $\calX$, from \eqref{est:eps:4}
we obtain that there exists a measure $\calT\in\calX'$ such that 
\begin{align}\label{limit:eps:5b}
  & \beta\ee(u^\epsi)\rightharpoonup \calT\quext{weakly star in }\,\calX'.
\end{align}
\end{subequations}
In view of the density of $\calX\cap \calV$ both in $\calX$ and in $\calV$, the
measure $\calT$ represents $\xi$ on $\calX$, i.e.~\eqref{identif} holds.

Let us now go back to \eqref{eqnweps}, now rewritten for general $t\in(0,T]$
and $\fhi\in \calV_t$. Then, rearranging terms, and using 
the above convergence relations~\eqref{limit:eps},
we obtain that there exists the limit 
\begin{align}\no
   \lim_{\epsi\searrow 0} 
    \duavw{ \beta\ee(u^\epsi),\fhi }\zt 
    & = \lim_{\epsi\searrow 0} \mbox{} \Big[  (\!( u^\epsi_t, \fhi_t )\!)\zt 
    -( u^\epsi_t(t) , \fhi(t) )
    - (\!( \nabla u^\epsi_t, \nabla \fhi )\!)\zt
    - (\!( \nabla u^\epsi, \nabla \fhi )\!)\zt\\
  \no 
   & \mbox{} ~~~~~~~~~~~~~~~
    + \lambda(\!( u^\epsi, \fhi )\!)\zt
    +( u_1\ee , \fhi(0) )
    + (\!( g, \fhi )\!)\zt \Big] \\
  \no
  & = (\!( u_t, \fhi_t )\!)\zt 
    - ( u_t(t) , \fhi(t) )
    - (\!( \nabla u_t, \nabla \fhi )\!)\zt
    - (\!( \nabla u, \nabla \fhi )\!)\zt\\[1mm]
  \label{eqnwepsx}  
   & \mbox{} ~~~~~~~~~~~~~~~
    + \lambda(\!( u, \fhi )\!)\zt
    + ( u_1 , \fhi(0) )
    + (\!( g, \fhi )\!)\zt.
\end{align}
A crucial point in our argument is that the left hand side tends, {\sl with no need
of extracting a further subsequence}, to a linear and continuous
functional on $\calV_t$ that acts on $\fhi$ as specified
by the \rhs. Noting as $\xi_{(t)}$ such a functional, we have
in other words
\begin{equation}\label{limit:eps:5c}
  \beta\ee(u^\epsi)\rightharpoonup \xi_{(t)}\quext{weakly in }\,\calV_t'.
\end{equation}
Moreover, \eqref{eqnwepsx} can be restated as
\begin{align}\no
   & \mbox{} - (\!( u_t, \fhi_t )\!)\zt 
    + ( u(t) , \fhi(t) )
    + (\!( \nabla u_t, \nabla \fhi )\!)\zt
    + (\!( \nabla u, \nabla \fhi )\!)\zt\\
  \label{eqnwlimit2}  
   & \mbox{} ~~~~~
    + \duavw{ \xi_{(t)},\fhi }\zt
    - \lambda(\!( u, \fhi )\!)\zt
    = ( u_1 , \fhi(0) )
    + (\!( g, \fhi )\!)\zt.
\end{align}
Hence, \eqref{eqnw} and \eqref{eqnwT}, which is a particular
case of it, are proved.
Note now that, from \eqref{est:eps:4}, it also follows
\begin{equation}\label{limit:eps:5c2}
  \beta\ee(u^\epsi) \rightharpoonup \calT_{(t)}\quext{weakly star in }\,\calX_t'
\end{equation}
and also this convergence holds with no need of extracting 
further subsequences. Indeed, the limit of the whole (sub)sequence
is already identified as $\xi_{(t)}$ on the dense subspace
$\calX_t\cap\calV_t$. This also implies that 
the measure $\calT_{(t)}$ represents $\xi_{(t)}$ on $\calX_t$ in the 
sense of \eqref{identif}. Using the fact that for any $\fhi\in\calV_{t,0}$
the extension $\tilfhi$ lies in $\calV$, it is easy to check
that the functionals $\xi$ and $\xi_{(t)}$ are ``compatible''.
Hence, we have checked points~{\sl (a)}\/ and~{\sl (b1)-(b2)}\/
of Definition~\ref{def:sol} of weak solution.

Let us now show relation \eqref{ener}, i.e., point~{\sl (e)}
of Definition~\ref{def:sol}. To this aim, we write
\eqref{eqnw} with $\fhi=u$ for $s,t\in(0,T]$ and take the difference.
Note that the choice $\fhi=u$ is admissible since $u\in\calV$.
We then infer
\begin{align}\no
   & \mbox{} 
   - \| u_t \|_{L^2(s,t;H)}^2
    + ( u_t(t) , u(t) )
    + \int_s^t ( \nabla u_t, \nabla u )\,\ditau
    + \| \nabla u \|_{L^2(s,t;H)}^2\\
  \label{ener0}
   & \mbox{} ~~~~~
    + \duavw{ \xi_{(t)},u }_{(0,t)}
    - \duavw{ \xi_{(s)},u }_{(0,s)}
    - \lambda \| u \|_{L^2(s,t;H)}^2
   = ( u_t(s) , u(s) )
    + \int_s^t \io g u\,\dix\,\ditau.
\end{align}
Then, computing explicitly the integral on the \lhs, \eqref{ener}
readily follows.

\medskip

\noindent%
\textbf{Step 4. Identification of $\xi$.}~~%
To conclude our proof we need to identify $\xi$ (and $\xi_{(t)}$) in the sense of the 
weak constraint~\eqref{incl2}. This will give~{\sl (c)}
of Definition~\ref{def:sol}. We start working on $\xi$, and,
to get the identification, we shall implement the so-called
Minty's trick in the duality between $\calV'$ and $\calV$. This corresponds
to checking the following two conditions:
\begin{itemize}
\item[(i)] There holds the $\limsup$-inequality
\begin{align}\label{claim}
  \limsup_{\epsilon\searrow 0} \,\duavw{\beta\ee(u^\epsi), u^\epsi} \leq \duavw{\xi, u}.
\end{align}
\item[(ii)] The operators $\beta\ee$ {\sl suitably}\/ converge to $\beta_w$,
  in such a way that~\eqref{incl2} may follow as a consequence of~\eqref{claim}.
\end{itemize}
\noindent%
We start by checking property~(i), postponing the discussion regarding
the correct notion of convergence for~(ii) and its implications. 
Writing \eqref{eqnweps} for $\fhi=u\ee$ and $t=T$, we obtain
\begin{align}\no
    \duavw{ \beta\ee(u^\epsi),u^\epsi }
     & = \| u^\epsi_t \|_{L^2(0,T;H)}^2
    - ( u^\epsi_t(T) , u^\epsi(T) ) 
    + ( u_1\ee , u_0\ee )
    - \frac12 \| \nabla u^\epsi(T) \|^2
    + \frac12 \| \nabla u_0\ee \|^2 \\
  \label{sci1}
   & \mbox{}~~~~~ 
   - \| \nabla  u^\epsi \|_{L^2(0,T;H)}^2
   + \lambda \| u^\epsi \|_{L^2(0,T;H)}^2
    + (\!( g, u\ee)\!).
\end{align}
Now, thanks to \eqref{limit:eps:1}, \eqref{limit:eps:2a}, \eqref{limit:eps:2b},
\eqref{limit:eps:3}, and \eqref{limit:eps:4}, we see that the $\limsup$ 
(as $\epsi\searrow 0$ along a proper subsequence)
of the right hand side is less or equal than
\begin{align}\no
   & \| u_t \|_{L^2(0,T;H)}^2
    - ( u_t(T) , u(T) ) 
    + ( u_1 , u_0 )
    - \frac12 \| \nabla u(T) \|^2
    + \frac12 \| \nabla u_0 \|^2
   - \| \nabla u \|_{L^2(0,T;H)}^2 \\
 \label{sci2}
  & \mbox{}~~~~~ 
   + \lambda \| u \|_{L^2(0,T;H)}^2
   + (\!(g, u)\!).
\end{align}
Hence, using \eqref{ener} written for $t=T$ and $s=0$, we see that
the above expression is equal to $\duavw{ \xi,u }$. Therefore,
\eqref{claim} is proved.

\smallskip

Let us now switch to discussing~(ii), which requires the introduction
of some additional machinery. We present it by following the lines
of the book by Attouch~\cite{At}. At first, we observe that 
the restriction to $\calV$ of the function $\beta\ee$ 
can be seen as a monotone operator from $\calV$ to $\calV'$
(once one works in the parabolic Hilbert triplet
$\calV \subset \calH \subset \calV'$). Indeed, if $v\in \calV$,
then $\beta\ee(v)\in \calH\subset\calV'$ by the 
Lipschitz continuity of $\beta\ee$. 
Hence, for any $v_1,v_2\in \calV$, we have
\begin{equation}\label{beta11}
  \duavw{\beta\ee(v_2)-\beta\ee(v_1), v_2- v_1}
   = (\!( \beta\ee(v_2)-\beta\ee(v_1), v_2- v_1 )\!)
   \ge 0.
\end{equation}
Moreover, if $v,z\in \calV$, then, by definition
of subdifferential,
\begin{equation}\label{beta12}
  \duavw{\beta\ee(v), z - v}
   = (\!( \beta\ee(v), z - v )\!)
   \le \iTo ( j\ee(z) - j\ee(v) )\,\dix\,\dit
   = \calJ\ee|_{\calV}(z) - \calJ\ee|_{\calV}(v).
\end{equation}
In other words, we have the graph inclusion
\begin{equation}\label{beta13}
  \beta\ee|_{\calV} 
   \subset \devv \calJ\ee|_{\calV},
\end{equation}
where the notation used on the \rhs\ stands for 
the subdifferential of $\calJ\ee|_{\calV}$ 
with respect to the duality
pairing between $\calV'$ and $\calV$. By the standard
theory of subdifferentials, this is a maximal monotone
operator from $\calV$ to $2^{\calV'}$, which includes
(in the sense of graphs) the (monotone, but not necessarily
maximal) operator $\beta\ee|_{\calV}$. 

In view of the fact that the family of functionals $\{\calJ\ee|_{\calV}\}$ 
(defined on $\calV$ and taking values in~$[0,+\infty)$) 
is increasing as $\epsi$ decreases
to $0$, applying \cite[Thm.~3.20]{At}, we obtain
\begin{equation}\label{beta14}
  \calJ\ee|_{\calV} \to \sup_{\epsi\in(0,1)} \calJ\ee|_{\calV}
\end{equation}
in the sense of {\sl Mosco convergence}\/ (that is Gamma-convergence
both in the strong and in the weak topology of~$\calV$). Moreover,
by the monotone convergence theorem it is readily seen that
the functional on the \rhs\ coincides in fact with $\calJ|_{\calV}$.
Hence, owing to \cite[Thm.~3.66]{At}, the family of maximal monotone
operators $\devv \calJ\ee|_{\calV}$, identified with the family
of their graphs in the product space $\calV\times \calV'$,
converges {\sl in the sense of graphs}\/ (cf.~\cite[Def.~3.58]{At})
to $\devv \calJ|_{\calV}=\beta_w$. Namely,
\begin{equation}\label{defGconv}
  \perogni [x;y] \in \devv \calJ|_{\calV}, \quad \esiste 
   [x\ee;y\ee] \in \devv \calJ\ee|_{\calV} \quext{such that \,
   $[x\ee;y\ee]\to[x;y]$~~strongly in }\, \calV\times \calV'.
\end{equation}
Hence, in view of the facts that $[u\ee;\beta\ee(u\ee)]\in \calJ\ee|_{\calV}$
(thanks to \eqref{beta13}), $[u\ee;\beta\ee(u\ee)] \to [u;\xi]$
weakly in $\calV\times \calV'$ (thanks to \eqref{limit:eps:1}
and \eqref{limit:eps:5}), and to the $\limsup$-inequality
\eqref{claim}, we may apply \cite[Prop.~3.59]{At}, yielding that
$[u;\xi] \in \devv \calJ\ee|_{\calV}=\beta_w$. Hence, \eqref{incl2}
is proved.

To conclude this part, we need to prove that 
$\xi_{(t)}\in \beta_{w,(t)}(u)$. To this aim, it is sufficient
to adapt the above argument by working on the 
subinterval~$(0,t)$. Indeed, relations~\eqref{eqnweps}
and \eqref{ener} (the latter for $s=0$) hold
on any subinterval $(0,t)$. Moreover, we can take advantage
of \eqref{limit:eps:1} (whose analogue
obviously holds also on subintervals) and \eqref{limit:eps:5c}.

\medskip

\noindent%
\textbf{Step 5. Further properties of solutions.}~~%
Let us start proving that inequality \eqref{claim} is
in fact an equality. Indeed, owing to \eqref{defGconv}, there exist
$[x\ee;y\ee] \in \devv \calJ\ee|_{\calV}$ such that 
$[x\ee;y\ee]\to[u;\xi]$ strongly in $\calV\times \calV'$.
Hence, noting that, by monotonicity,
\begin{equation}\label{claim2}
  0 \le \duavw{\beta\ee(u^\epsi)-y\ee, u^\epsi-x\ee},
\end{equation}
taking the $\liminf$ as $\epsi\searrow 0$, and recalling
\eqref{claim}, we obtain
\begin{equation}\label{claim3}
  \lim_{\epsilon\searrow 0} \,\duavw{\beta\ee(u^\epsi), u^\epsi} 
   = \duavw{\xi, u}.
\end{equation}
As a consequence, the limit of the \rhs\ of \eqref{sci1} exists and
coincides with \eqref{sci2}. In view of the fact that convergence of 
most terms of \eqref{sci1} is already known from the previous estimates
we get in particular that
\begin{equation}\label{claim4}
 \lim_{\epsi\searrow0} \bigg( \frac12\| \nabla u^\epsi(T) \|^2  
      + \| \nabla  u^\epsi \|_{L^2(0,T;H)}^2 \bigg)
  = \frac12\| \nabla u(T) \|^2 + \| \nabla  u \|_{L^2(0,T;H)}^2.
\end{equation}
As before, this argument can be repeated on any subinterval $(0,t)$. 
Hence, recalling \eqref{limit:eps:1} and \eqref{limit:eps:2b}, 
we finally arrive at
\begin{align}\label{limit:eps:7}
  & u^\epsi\rightarrow u \quext{strongly in }\,L^2(0,T;V),\\
 \label{limit:eps:7*}
  & u^\epsi(t) \rightarrow u(t) \quext{strongly in }\, V
   ~~\text{for all }\,t\in[0,T].
\end{align}

\smallskip

Next, let us show that, under assumption \eqref{hp:u0:2},
the additional regularity \eqref{regou2} holds. To this aim,
we go back to the approximate problem, and, 
in the spirit of~\cite{PZ}, we test \eqref{eqne} by $Au\ee$. 
Indeed, $u\ee$ has sufficient smoothness in order for this procedure
to be admissible (cf.~\eqref{regouee}).
Integrating by parts, and using the monotonicity of $\beta\ee$, 
we then easily infer
\begin{equation}\label{patzel1}
  \ddt ( u_t\ee , A u\ee )
   + \frac12 \ddt \| Au\ee \|^2 
   + \| Au\ee \|^2
   = ( \lambda u\ee + g, Au\ee )
   + ( u_t\ee , A u\ee_t ).
\end{equation}
By some further integration by parts and using H\"older's and Young's
inequalities (and the definition of the operator $A$), the \rhs\ can be easily estimated as follows:
\begin{equation}\label{patzel1b}
  ( \lambda u\ee + g, Au\ee )
   + ( u_t\ee , A u\ee_t )
   \le \frac12 \| Au\ee \|^2 
   + C \big( \| g \|^2 + \| u\ee \|^2 \big)
   + \| \nabla u\ee_t \|^2. 
\end{equation}
Here and below, $C>0$ is a constant independent of $\epsi$.
Hence, integrating \eqref{patzel1} over $(0,t)$ for arbitrary 
$t\in(0,T]$, using~\eqref{patzel1b},
and recalling estimate~\eqref{est:eps:1}, we
easily obtain
\begin{equation} \label{patzel4}
  2 (u\ee_{t}(t), Au\ee(t))
   + \| Au\ee(t) \|^2 
   + \itt \| Au\ee \|^2 \,\dis
 \le C + \| A u_0\ee \|^2.
\end{equation}
Now, one can immediately check that, 
under assumption \eqref{hp:u0:2}, if 
$u_0\ee$ is defined as in \eqref{ellee},
then the \rhs\ of \eqref{patzel4} is bounded independently 
of $\epsi$. Hence, noticing that the \lhs\ is larger or
equal than
\begin{equation} \label{patzel5}
  \frac12 \| Au\ee(t) \|^2 
   - C \| u\ee_t(t) \|^2 
   + \itt \| Au\ee \|^2 \,\dis,
\end{equation}
where the second term is uniformly controlled due 
to~\eqref{est:eps:1}, we readily arrive
at 
\begin{equation} \label{patzel6}
 \| u\ee \|_{L^\infty(0,T;D(A))} \le M.
\end{equation}
Letting $\epsi\searrow0$, we then infer
\begin{equation} \label{patzel6b}
  u\in L^\infty(0,T;D(A))
\end{equation}
thanks to semicontinuity of norms with respect to weak convergence. 
Finally, \eqref{regou2}, i.e., weak continuity of $u$ with values
in $D(A)$, follows by combining \eqref{patzel6b} with the regularity 
$u\in C([0,T];V)$ following from~\eqref{regou}, 
and applying standard results. 

\medskip

Eventually, we show that weak solutions constructed
as limit points of $\{u\ee\}$ also satisfy a form of the energy
{\sl inequality}. We start by proving it on intervals
of the form $[0,t]$, $t\in(0,T]$.
To this aim, we write relation \eqref{enereps} for $s=0$ and 
take the $\liminf$ as $\epsi\searrow0$. 
Then, using standard
semicontinuity arguments together with relations
\eqref{inizee1}-\eqref{inizee2}, \eqref{limit:eps:1},
and \eqref{limit:eps:2a}-\eqref{limit:eps:2b}, it is not difficult
to infer, for every $t\in(0,T]$,
\begin{align}\no
  & \mbox{} 
   \frac{1}{2}\| u_t(t) \|^2
   + \int_\Omega j(u(t)) \,\dix
   - \frac{\lambda}{2} \|u(t)\|^2
   + \frac{1}{2}\| \nabla u (t)\|^2
   + \| \nabla u_t \|_{L^2(0,t;H)}^2\\
 \label{enerlim}
  & \mbox{} ~~~~~
   \le \frac{1}{2}\| u_1 \|^2
   + \int_\Omega j(u_0)\,\dix
   - \frac{\lambda}{2} \|u_0 \|^2
   + \frac{1}{2}\| \nabla u_0 \|^2
   + \int_0^t ( g,u_t )\,\ditau.
\end{align}
Note in particular that relation 
\begin{equation}\label{enerlim2a}
  J(u(t)) = \int_\Omega j(u(t)) \,\dix
   \le \liminf_{\epsi\searrow0}
   \int_\Omega j\ee(u\ee(t)) \,\dix
   = \liminf_{\epsi\searrow0}
   J\ee(u\ee(t)) 
\end{equation}
is a consequence of \eqref{limit:eps:2a} and of the fact 
that the functionals $J\ee$ converge to $J$ in the 
sense of Mosco (cf.~\cite[Par.~3.3]{At})
in the space~$H$. Recalling \eqref{defiE},
\eqref{enerlim} reduces to \eqref{enerlim2} in the
case $s=0$. Let us now consider a generic
interval $[s,t]$ for $0 < s < t \le T$ and let 
us go back to~\eqref{enereps} written for this choice of~$s,t$. 
Let us take once more the $\liminf$ as $\epsi\searrow 0$. 
Then, the \lhs\ is treated exactly as before.
On the other hand, when looking at the \rhs,
it is easy to check that
\begin{equation}\label{ener21}
    - \frac{\lambda}{2} \|u^\epsi(s)\|^2
    + \frac{1}{2}\| \nabla u^\epsi (s)\|^2
    + \int_s^t ( g,u^\epsi_t )\,\ditau
  \rightarrow  
    - \frac{\lambda}{2} \|u(s)\|^2
    + \frac{1}{2}\| \nabla u(s) \|^2
    + \int_s^t ( g,u_t )\,\ditau 
\end{equation}
thanks in particular to \eqref{limit:eps:1},
and \eqref{limit:eps:7*}.
Next, thanks to \eqref{limit:eps:4}, and up to 
extracting a further subsequence of $\epsi\searrow0$,
we have
\begin{equation}\label{ener22}
  \frac{1}{2}\| u^\epsi_t(s) \|^2
   \to \frac{1}{2}\| u_t(s) \|^2
\end{equation}
for {\sl almost}\/ every choice of $s\in(0,T)$. Next, we need to 
control the component of the energy related with the constraint.
Namely, we would like to prove that, at least for
a.e.~$s\in(0,T)$,
\begin{equation}\label{ener23}
  \limsup_{\epsi\searrow0} \int_\Omega j\ee(u^\epsi(s))\,\dix
   \le \int_\Omega j(u(s))\,\dix
\end{equation}
(hence, coupling the above with \eqref{enerlim2a} written for $t=s$, 
we would get convergence of that term). We can start noticing that 
\begin{equation}\label{ener24}
  \int_\Omega j\ee(u^\epsi(s))\,\dix
   = \int_\Omega \big( j\ee(u\ee(s)) - j\ee(u(s)) \big)\,\dix
   + \int_\Omega j\ee(u(s))\,\dix.
\end{equation}
Moreover,
\begin{equation}\label{ener25}
  \int_\Omega j\ee(u(s))\,\dix
   \to \int_\Omega j(u(s))\,\dix
\end{equation}
by the monotone convergence theorem. Now, by definition of 
subdifferential, we may write
\begin{equation}\label{ener26}
  \int_\Omega \big( j\ee(u\ee(s)) - j\ee(u(s)) \big)\,\dix
   \le \int_\Omega \beta\ee(u\ee(s)) ( u\ee(s) - u(s) ) \,\dix
   =: \mu\ee(s),
\end{equation}
and we have to discuss the behavior of the functions~$\mu\ee$.
First, we observe that
\begin{equation}\label{ener27}
  \int_0^T \mu\ee(s) \,\dis
  = \duavw{ \beta\ee(u\ee), u\ee - u }  
  \to 0,
\end{equation}
the latter property following from \eqref{claim3} and \eqref{limit:eps:5}.
Moreover, thanks to \eqref{ener26}, we have 
\begin{equation}\label{ener26b}
  \mu\ee(s) \ge J\ee(u\ee(s)) - J\ee(u(s))
    \ge J\ee(u\ee(s)) - J(u(s)),
\end{equation}
whence $\liminf_{\epsi\searrow0} \mu\ee(s) \ge 0$ thanks 
to \eqref{enerlim2a}. Hence, we have in particular
$\lim_{\epsi\searrow0} (\mu\ee)^-(s) = 0$, 
$(\cdot)^-$ denoting the negative part.  Moreover,
from \eqref{ener26b}, \eqref{enerlim2a} and \eqref{est:eps:3} we infer
\begin{equation}\label{ener26c}
  \mu\ee(s) \ge - J(u(s))
   \ge - \liminf_{\delta\searrow0} J^\delta(u^\delta(s)) 
   \ge - M,
\end{equation}
for all $\epsi\in(0,1)$ and $s\in(0,T]$. Hence, by the dominated convergence
theorem we obtain that $(\mu\ee)^-\to 0$ in $L^1(0,T)$. 
Consequently, thanks to \eqref{ener27}, we conclude
that $\mu\ee\to 0$ in $L^1(0,T)$. 
Hence, up to a subsequence, $\mu\ee \to 0$ almost everywhere in $(0,T)$,
whence \eqref{ener23} follows. This actually implies
\eqref{enerlim2} for {\sl almost}\/ every $s\in(0,T)$
and every $t\in(s,T]$, as desired. The proof of
Theorem~\ref{teo:esi} is concluded.
%
%
\beos\label{rem:dissi}
 If the source term $g$ is $0$,
 from \eqref{enerlim2} follows in particular that 
 the energy loss in the time interval $(s,t)$ is
 {\sl at least}\/ as large as the dissipation term
 $\calD(s,t):= \| \nabla u_t \|_{L^2(s,t;H)}^2$. 
 Of course, as commonly occurs situations characterized
 by bad regularity, the energy dissipated 
 may be in fact strictly larger than $\calD(s,t)$.
 Indeed, we may observe that
 proving equality in \eqref{enerlim} 
 appears out of reach in the present regularity
 setting. 
\eddos
\beos\label{rem:reg2}
 In view of our strategy of proof for Theorem~\ref{teo:esi}, we can 
 give some further observation complementing Remark~\ref{rem:reg1}. 
 Hence, let us go back to the ``toy problem'' \eqref{toy},
 for example with $\beta = \de I_{[-1,1]}$
 (but our consideration also apply to different choices of $\beta$).
 Then, implementing our regularization method we get the
 equation
 \begin{equation}\label{ODEee}
   u\ee_{tt} + \beta\ee(u\ee) = 0.
 \end{equation}
 Setting $v\ee:=u\ee_t$, the shape of solution trajectories
 of the 2D ODE system associated to \eqref{ODEee} in the phase
 space $(u\ee,v\ee)$ can be easily described. In particular, 
 since (in this spatially homogeneous setting) no dissipation occurs,
 trajectories are periodic. Moreover, we may notice that,
 for $\epsi\searrow0$, $(u\ee,v\ee)$ converges in a suitable way
 to a couple $(u,v)$, where $v=u_t$ and $u$ 
 solves \eqref{toy}. Clearly, $(u,v)$ is also 
 a periodic trajectory and its image in the phase space
 lies in some level set $\{j(u)+v^2/2=c\}$, $c\ge 0$, 
 of the ``energy'' functional. In particular, 
 whenever $D(j)=[-1,1]$ (as happens in the case of the 
 indicator function $j=I_{[-1,1]}$, and also for the ``logarithmic
 potential'' mentioned in the introduction),
 such level sets are (at least for large initial energy,
 i.e., for large values of $c$)
 not connected. Namely, their shape determines the 
 jumps of~$u_t$ (which, consequently, cannot
 occur in an ``arbitrary'' way). Note also that taking
 different choices for the approximations $\beta\ee$ of $\beta$ 
 does not modify the shape of $(u,v)$. Of course it is clear
 that, in the case of our equation \eqref{eqn:intro},
 the situation is much more complicated than for
 \eqref{toy} in view of the infinite-dimensional setting. 
 However, the fact that our weak
 solutions $u$ are still
 built as limit points of families
 $u\ee$ solving a very natural regularization
 of the equation suggests that the jumps of 
 $u_t$ occurring in the limit may be in some sense 
 ``physical'', i.e., they are determined by the fact
 that $u\ee_t$, as $\epsi\searrow0$, may tend
 to develop discontinuities. In other words, the occurrence of
 ``spurious'' jumps of $u_t$ (as are the somehow 
 ``arbitrary'' jumps described in Remark~\ref{rem:reg1})
 should be excluded in view of the fact that
 our weak solutions descend from the approximation scheme.
\eddos
\beos\label{rem:tempo2}
 Let us give some further observation complementing
 Remark~\ref{rem:tempo}. Again, we consider, just 
 for simplicity, the ``toy'' model~\eqref{toy}; 
 however, our considerations also apply to
 the original equation~\eqref{eqn:intro}.
 Actually, from our approximation
 argument we know that, for any $t\in (0,T]$,
 (a subsequence of) $\beta\ee(u\ee)$ (weakly star) 
 converges to a measure $\calT_{(t)}$
 on $\ov{Q_t}$ (in particular, we have convergence
 to some $\calT$ on the whole interval). 
 In the toy case, of course, $\ov{Q_t}=[0,t]$; moreover,  
 we are allowed to identify 
 $\calT_{(t)}=\xi_{(t)}\in \calV_t'$
 because Sobolev functions are continuous in~1D.
 Let us now consider the particular case when
 $\beta\ee(u\ee)$ is supported in
 some interval of the form $[t\ee-\epsi,t\ee+\epsi]$
 and is $0$ outside that interval. Then, assuming that
 $t\ee$ converges to some point $t\in(0,T)$ as $\epsi\searrow 0$,
 and $\beta\ee(u\ee)$ ``spikes'' around $t\ee$ in a proper way,
 it may happen that $\beta\ee(u\ee)$ (weakly star) converges
 to $\calT=\delta_t$ (the Dirac
 delta concentrated in $t$) in $\calX'$.
 This kind of behavior may be (possibly) driven for instance 
 by inserting a nonzero forcing term~$g$ in the equation.
 Then, in the case when, for instance, $t\ee = t-2\epsi$,
 it turns out that the singularity of $\calT$ develops 
 {\sl before}~$t$. Consequently, 
 $\beta\ee(u\ee)$ also converges to $\delta_t$ 
 in $\calX_t'$. In particular, \eqref{eqnw}
 holds in $[0,t]$ with $\xi_{(t)}=\delta_t$. On the other hand,
 if $t\ee = t+2\epsi$, i.e.,
 the singularity of $\calT$ develops 
 {\sl after}~$t$, in that case
 $\beta\ee(u\ee)$ converges to $\calT_{(t)}=0$ 
 in $\calX_t'$, whence \eqref{eqnw}
 holds in $[0,t]$ with $\xi_{(t)}=0$. Note that
 this happens in spite of the fact that the limit 
 measure~$\calT$ over the whole $[0,T]$ is the same in the 
 two cases. This fact suggests that
 the formulation \eqref{eqnw} on the subinterval~$(0,t)$
 contains some additional information that cannot be simply
 inferred by restricting the {\sl global}~formulation
 \eqref{eqnwT}. This is the reason why we decided 
 to include~{\sl (b2)}\/ in our existence theorem.
\eddos
\beos\label{rem:reg3}
 Let $t$ be one of the (at most countably many) jump points of~$u_t$.
 Then, both the point value of~$u_t$ at $t$ and the occurrence
 of concentration phenomena for the measure $\calT_{(t)}$ at the same point
 also depend on the choice of the approximating
 problem (i.e., of $\beta\ee$; actually our argument
 works provided that $\beta\ee$ is smooth and converges to $\beta$ in
 the graph sense) and of the selection of converging
 subsequences via Helly's theorem. This can be seen again by looking
 at the ``toy equation''~\eqref{toy} with $\beta=\de I_{[-1,1]}$
 and initial values $u_0=0$ and $u_1=1$.
 Then, we know that the (first) jump of $u_t$ occurs at $t=1$.
 Let us now consider the approximation
 \eqref{ODEee} with the choice
 \begin{equation}\label{betaee1}
   \beta\ee(r)= \begin{cases}
     0 & \text{~~if }|r|\le r\ee, \\
     \epsi^{-2}(r-r\ee) & \text{~~if } r > r\ee, \\
     \epsi^{-2}(r+r\ee) & \text{~~if } r < - r\ee, 
    \end{cases}
 \end{equation}
 where, for any $\epsi\in(0,1)$, one may choose 
 (in an arbitrary way) $r\ee$ in the interval
 $[1-\epsi\pi,1]$. It is then clear that, whatever 
 are the chosen values of $r\ee$, $\beta\ee$ tends
 to $\beta=\de I_{[-1,1]}$ in the sense of graphs 
 as $\epsi\searrow0$. Hence, our limit problem is 
 the desired one. Let us notice that, for $r\ge 0$,
 we have
 \begin{equation}\label{betaee2}
   j\ee(r) = \frac1{2\epsi^2} \big((r-r\ee)^+\big)^2.
 \end{equation}
 In particular, $j\ee(1)=0$ if $r\ee=1$, whereas 
 \begin{equation}\label{betaee3}
   j\ee(1)=\frac{(1-r\ee)^2}{2\epsi^2}
    \in \Big(0,\frac{\pi^2}2\Big] \quext{if }\,r\ee\in[1-\epsi\pi,1).
 \end{equation}
 Then, $u\ee(1)=1$ if $r\ee=1$,
 whereas in case $r\ee\in[1-\epsi\pi,1)$ one can easily compute
 $u\ee(t) = r\ee + \epsi \sin \left(\frac{t-r\ee}{\epsi}\right)$
 for $t\in (r\ee,r\ee+\epsi\pi]$, whence 
 $u\ee_t(t)= \cos \left(\frac{t-r\ee}{\epsi}\right)$
 and $u\ee_t(1)= \cos \left(\frac{1-r\ee}{\epsi}\right)$.
 Hence, choosing appropriately (and somehow ``wildly'') 
 $r\ee$ in the interval $[1-\epsi\pi,1]$
 as $\epsi$ varies in $(0,1)$,
 one may obtain the effect that for any number
 $\ell\in[-1,1]$ there exists a subsequence $\epsi_n\searrow 0$ such that 
 $u^{\epsi_n}_t(1)$ tends to~$\ell$. 
 The use of Helly's theorem selects {\sl one}\/ of 
 these subsequences and determines the limit value $u_t(1)=\ell$
 (and, in turn, how the limit measures $\calT_{(s)}$
 concentrate at the jump point $t=1$).
\eddos
\beos\label{tempo3}
 As observed in the previous Remark, $\xi_{(t)}$ is not represented, in general,
 by the restricion of the measure $\calT$ to the set $\ov{Q_{t}}$. However, we can 
 give a more explicit characterization of this restriction in the following sense. From 
 \eqref{limit:eps:2c*} we have $u_t\in BV(0,T;X)$,
 where we may assume $X$ be the dual of a separable space.
 Hence, for all times $t\in[0,T]$ there exists (in the weak star topology
 of $X$) the limit
 \begin{align}
   u_t(t^+):=w^*\!\!\!-\!\!\!\lim_{s\rightarrow t^+}u_t(s).
 \end{align}
 Moreover this value coincides with the weak star limit
 \begin{align}
   w^*\!\!\!-\!\!\!\lim_{s\rightarrow t^+}\frac{1}{s-t}\int_t^s u_t(r)\,\dir.
 \end{align}
 In particular the limits above must hold with respect to the weak topology 
 of~$H$, since $u_t$ is bounded in $H$ uniformly in time. Let us now write \eqref{eqnwT} 
 with $\varphi$ replaced by $\varphi h_s\in \mathcal X$, with $\varphi\in\mathcal X\cap \calV$
 and $h_s:[0,T]\rightarrow [0,1]$ be the function such that $h_s=1$ on $[0,t]$, 
 $h_s=0$ on $[s,T]$, and $h_s$ be affine in $[t,s]$. We obtain
 \begin{align}\no
   & \mbox{} - (\!( u_t, \fhi_t )\!)_{(0,t)}
   -(\!( u_t, \fhi_t h_s)\!)_{(t,s)}
   +\frac{1}{s-t}(\!( u_t, \fhi)\!)_{(t,s)} 
    + (\!( \nabla u_t, \nabla (\fhi h_s) )\!) 
    + (\!( \nabla u, \nabla (\fhi h_s) )\!)\\
  \no  
   & \mbox{} ~~~~~
    + \iint_{\ov{Q}} \fhi h_s \,\deriv\!\calT
   - \lambda (\!( u, \fhi h_s)\!) 
    = ( u_1 , \fhi(0) )
    + (\!( g, \fhi h_s)\!).
 \end{align}
 Letting $s\searrow t$, we see that the third term tends to $( u_t(t^+) , \fhi(t) )$, 
 while the other terms pass to the limit thanks to the dominated convergence theorem 
 and the fact that $\varphi h_s\rightarrow\varphi\chi_{[0,t]}$ pointwise, 
 so in particular $\mathcal T$-almost everywhere. We then obtain
 \begin{align}\no
   & \mbox{} - (\!( u_t, \fhi_t )\!)_{(0,t)}+( u_t(t^+), \fhi(t)) 
    + (\!( \nabla u_t, \nabla \fhi)\!)_{(0,t)} 
    + (\!( \nabla u, \nabla \fhi)\!)_{(0,t)}\\
  \no  
   & \mbox{} ~~~~~
   + \iint_{\ov{Q}} \fhi \,\deriv( \calT \chi_{[0,t]})
    - \lambda (\!( u, \fhi)\!)_{(0,t)} 
    = ( u_1 , \fhi(0) )
    + (\!( g, \fhi)\!)_{(0,t)}.
 \end{align}
 Comparing with \eqref{eqnw}, we deduce that $\xi_{(t)}$ is represented 
 by the restriction of $\mathcal T$ to the closed set $\ov{Q_t}$
 whenever the pointwise value $u_t(t)$ coincides with 
 $u_t(t^+)$, which happens in fact in the complementary of a countable set of times.
 In other words, in that case we have $\mathcal T\llcorner_{\ov{Q_t}}=\mathcal T_{(t)}$.
 %
 %
 %
\eddos
\beos\label{energia2}
 Relation \eqref{enerlim2} implies in particular that, at least when
 $g\equiv 0$, the energy functional coincides almost everywhere with a 
 nonincreasing function. In a sense this fact provides an additional 
 criterion for selecting which are the ``admissible'' jumps of $u_t$
 (cf.~Remark~\ref{rem:reg2}). Namely, jumps may
 occur only in such a way that they do not increase 
 the total energy of the system. For $g\not=0$ similar considerations
 hold, up to the fact that $g$ acts somehow as an additional energy
 source. 
\eddos
\beos\label{energia3}
 It is maybe also worth stressing that Theorem~\ref{teo:esi}
 states the existence of {\sl at least one}\/ weak solution satisfying the 
 properties detailed above. Due to nonuniqueness, there may well 
 exist ``spurious'' solutions having worse properties. For example
 they may be constructed in such a way that the time derivative 
 $u_t$ admits somehow ``nonphysical'' jumps. However our procedure
 shows that {\sl every}\/ weak solution that is a limit point
 of our natural regularization scheme is ``physical''
 (for example, in view of \eqref{enerlim2}, energy-increasing 
 jumps cannot occur).
\eddos


\section*{Acknowledgements}

The financial support of the FP7-IDEAS-ERC-StG \#256872
(EntroPhase) is gratefully acknowledged by the authors. The present paper 
also benefits from the support of the MIUR-PRIN Grant 2010A2TFX2 ``Calculus of Variations''
for EB and GS, and the GNAMPA (Gruppo Nazionale per l'Analisi Matematica, la Probabilit\`a 
e le loro Applicazioni) of INdAM (Istituto Nazionale di Alta Matematica). 
%



\end{document}